\documentclass{article}
\usepackage{amsmath}
\usepackage{amssymb}
\usepackage{graphicx,psfrag}
\newenvironment{proof}{\noindent {\bf Proof}}
{\hfill $\square$ \vspace{0.25cm}}

\newcommand{\ds}{\displaystyle}
\newcommand{\intot}{\ds\int_0^t}
\newcommand{\rr}{{\mathbb{R}}}
\newcommand{\nn}{{\mathbb{N}}}
\newcommand{\ala}{\nonumber \\}

\newcommand{\indiq}{1\!\! 1}
\newcommand{\e}{{\varepsilon}}
\newcommand{\mes}{{\cal M}_f^+}
\newcommand{\lc}{\left<}
\newcommand{\rc}{\right>}

\newcommand{\xhi}{{{\cal X}}}
\newcommand{\nk}{{n_k}}
\newcommand{\smn}{\mathcal{S}(n,m)}

\newtheorem{theo}{\indent Theorem}[section]
\newtheorem{prop}[theo]{\indent Proposition}

\newtheorem{lemma}[theo]{\indent Lemma}
\newtheorem{defin}[theo]{\indent Definition}
\newtheorem{cor}[theo]{\indent Corollary}

\begin{document}

\title{Marcus-Lushnikov processes, Smoluchowski's and 
Flory's models}

\author{Nicolas {\sc Fournier}\footnote{
Centre de Math\'ematiques,
Facult\'e de Sciences et Technologie,
Universit\'e Paris~XII, 61 avenue du G\'en\'eral de Gaulle, F--94010 Cr\'eteil 
Cedex, France. E-mail: {\tt nicolas.fournier@univ-paris12.fr}} 
\  and Philippe {\sc Lauren\c cot}\footnote{Institut de
Math\'ematiques de Toulouse, CNRS UMR~5219, Universit\'e Paul Sabatier
(Toulouse~III), 118 route de Narbonne, F--31062 Toulouse cedex 9,
France. E-mail: {\tt laurenco@mip.ups-tlse.fr}} } 

\date{}

\maketitle

\def\abstractname{Abstract}
\begin{abstract}
The Marcus-Lushnikov process is a finite stochastic particle system
in which each particle is entirely characterized by its mass.
Each pair of particles with masses $x$ and $y$ merges into
a single particle at a given rate $K(x,y)$.
We consider a {\it strongly gelling} kernel behaving as $K(x,y)=x^\alpha y
+ x y^\alpha$ for some $\alpha\in (0,1]$. In such a case, 
it is well-known that {\it gelation} occurs, that is, giant particles emerge.
Then two
possible models for hydrodynamic limits of the Marcus-Lushnikov process arise:
the Smoluchowski equation, in which the giant particles are
inert, and the Flory equation, in which the giant particles
interact with finite ones.\\
We show that, when using a suitable cut-off coagulation kernel in the
Marcus-Lushnikov process and letting the number of particles increase
to infinity, the possible limits solve either the Smoluchowski equation
or the Flory equation.\\
We also study the asymptotic behaviour of the largest particle in the 
Marcus-Lushnikov process without cut-off and show that there is only
one giant particle. This single giant particle represents,
asymptotically, the lost mass of the solution to the Flory equation.
\end{abstract}


\bigskip



\section{Introduction}\label{intro} \setcounter{equation}{0}

We investigate the connection between a stochastic coalescence model,
the Marcus-Lushnikov process, and two deterministic coagulation
equations, the Smoluchowski and Flory equations. Recall that the
Marcus-Lushnikov process \cite{lushnikov,marcus} is a finite
stochastic system of coalescing particles while the Smoluchowski and
Flory equations describe the evolution of the concentration $c(t,x)$
of particles of mass $x\in (0,\infty)$ at time $t\geq 0$ in an infinite
system of coalescing particles. Both models depend on a {\it
coagulation kernel} $K(x,y)$ describing the likeliness that two
particles with respective masses $x$ and $y$ coalesce. When $K$
increases sufficiently rapidly for large values of $x$ and $y$, a
singular phenomenon known as {\it gelation} occurs: giant particles
(that is, particles with infinite mass) appear in finite time (see
Jeon \cite{jeon}, Escobedo-Mischler-Perthame \cite{emp}). There is
however a strong difference between the Smoluchowski and Flory
equations: for the former, the giant particles are inert, while for
the latter, the giant particles interact with the finite particles.

When $K(x,y)/y\longrightarrow 0$ as $y\to\infty$ for all $x\in (0,\infty)$,
it is by now well-known that the Marcus-Lushnikov process converges to
the solution of the Smoluchowski equation when the number of particles
increases to infinity (see, e.g., Jeon \cite{jeon} and Norris
\cite{norris}). On the other hand, it has been shown in \cite{fg}
that, if $K(x,y)/y\longrightarrow l(x)\in (0,\infty)$ as $y\to\infty$
for all $x\in (0,\infty)$, then the Marcus-Lushnikov process converges
to the solution of the Flory equation.

Our aim in this paper is to study more precisely how this transition
from the Smoluchowski equation to the Flory equation arises in the
Marcus-Lushnikov process. For a coagulation kernel $K$ of the form
$K(x,y) \simeq x y^\alpha + x^\alpha y$ for some $\alpha \in (0,1]$,
we consider a Marcus-Lushnikov process starting with $n$ particles, 
with total mass $m_n$, where coalescence between particles larger than
some threshold mass $a_n$ is not allowed. We show that, in the limit
of large $n$, $m_n$ and $a_n$, this Marcus-Lushnikov process converges, 
up to extraction of a subsequence, either to the solution of the Flory equation
or that of the Smoluchowski equation, according to the behaviour of
$a_n/m_n$ for large values of $n$.

We also study the behaviour of the largest particles in the 
Marcus-Lushnikov process without cut-off, and show that, in some sense,
the total lost mass of the Flory equation is represented by one
giant particle in the Marcus-Lushnikov process. 
Aldous \cite{aldousgp} proved other results about giant particles 
for some similar (but more restrictive) kernels. We in fact obtain a
much more precise result about the size of the largest particle after
gelation, but we are not able to extend to our class of kernels his
result about the largest particle before gelation.


\section{Main result}\label{resultat}\setcounter{equation}{0}

Throughout the paper, a {\it coagulation kernel} is a function
$K: (0,\infty)^2 \mapsto [0,\infty)$ such that $K(x,y)=K(y,x)$
for all $(x,y) \in (0,\infty)^2$. 
We denote by $\mes$ the set
of non-negative finite measures on $(0,\infty)$.
Let us first recall the definition of the Marcus-Lushnikov process.

\begin{defin}\label{dfml}
Consider a coagulation kernel $K$, and an initial state 
$\mu_0 = m^{-1} \sum_{i=1}^n \delta_{x_i}$, with $n\geq 1$, 
$(x_1,...,x_n)\in(0,\infty)^n$ and $m=x_1+...+x_n$.
A c\`adl\`ag $\mes$-valued Markov process $(\mu_t)_{t\geq 0}$ is a 
Marcus-Lushnikov process associated with the 
pair $(K,\mu_0)$ if it a.s. takes its values in
\begin{equation}
\smn := \left\{ \frac{1}{m}\ \sum_{i=1}^k \delta_{y_i}, \; 1\le 
k\leq n, \;  (y_i)_{1\le i\le k} \in (0,\infty)^k, 
\sum_{i=1}^k y_i = m  \right\}
\end{equation}
and its generator is given by
\begin{eqnarray}\label{gml}
L^{K,\mu_0} \psi (\mu) =
\sum_{i\ne j}\left\{
\psi\left[ \mu + m^{-1}\left(\delta_{y_i+y_j} 
- \delta_{y_i} -\delta_{y_j} \right)\right] 
- \psi\left[\mu \right] \right\}
\frac{K(y_i,y_j)}{2m}
\end{eqnarray}
for all measurable functions $\psi:\mes\mapsto \rr$ and all 
states $\mu=m^{-1} \sum_{i=1}^k \delta_{y_i} \in \smn$.
\end{defin}

This process is known to be well-defined and unique, 
without any assumption on $K$, see, e.g., Aldous
\cite[Section 4]{aldous} or Norris \cite[Section 4]{norris}. 

\medskip

We now describe the Smoluchowski and Flory coagulation equations and
first introduce the class of coagulation kernels to be considered in
the sequel. As already mentioned, we will deal with kernels of the
form $K(x,y)\simeq x^\alpha y + x y^\alpha$ for some $\alpha \in
(0,1]$. More precisely, we assume the following:

\medskip

\noindent\textbf{Assumption $(A_\alpha)$}: The coagulation kernel $K$
is continuous on $(0,\infty)^2$ and there are $\alpha \in (0,1]$,
$l\in\mathcal{C}((0,\infty))$, and positive real numbers
$0<c<C<\infty$ such that 
$$
\lim_{y \to \infty} {K(x,y)}/{y} = l(x)\,,
$$
and 
\begin{equation}\label{assumption}
c\ (x^\alpha y + x y^\alpha)\leq K(x,y) \leq C\ (x^\alpha y + x
y^\alpha)  \;\;\mbox{ and }\;\;  c\ x^\alpha \leq l(x) \leq C\
x^\alpha 
\end{equation}
for all $(x,y) \in (0,\infty)^2$. 

\medskip

For such coagulation kernels, weak solutions to the Smoluchowski and
Flory coagulation equations are then defined as follows: 

\begin{defin}
Consider a coagulation kernel $K$ satisfying $(A_\alpha)$ for some
$\alpha \in (0,1]$ and  
$\mu_0\in \mes$ such that $\left< \mu_0(dx), 1+x\right>
< \infty$. For $\phi:(0,\infty) \mapsto \rr$, set
\begin{equation}\label{dfdelta}
\Delta \phi(x,y):= \phi(x+y)-\phi(x)-\phi(y)\,.
\end{equation}
A family $(\mu_t)_{t\geq 0} \subset \mes$ such that
$t \mapsto \left< \mu_t(dx), x\right>$ and $t\mapsto \lc \mu_t(dx),1\rc$
are non-increasing is a solution to:\\
(i) the Smoluchowski equation $(S)$ if
\begin{eqnarray}\label{smolu}
\left< \mu_t, \phi \right> &=& \left< \mu_0, \phi \right>
+ \frac{1}{2} \intot \left< \mu_s(dx) \mu_s(dy),  K(x,y)\Delta \phi(x,y) 
\right> ds
\end{eqnarray}
for all $\phi \in C_c([0,\infty))$ and $t\geq 0$;\\
(ii) the Flory equation $(F)$ if for 
\begin{eqnarray}\label{flory}
\left< \mu_t, \phi \right> &=& \left< \mu_0, \phi \right>
+ \frac{1}{2} \intot \left< \mu_s(dx) \mu_s(dy), K(x,y)\Delta\phi(x,y)
\right> ds
\ala
&&- \intot  \left< \mu_s(dx), \phi(x) l(x) \right> \left<\mu_0(dx)-\mu_s(dx), 
x \right> ds
\end{eqnarray}
for all $\phi \in C_c([0,\infty))$ and $t\geq 0$. Here and below,
$C_c([0,\infty))$ denotes the space of continuous functions with
compact support in $[0,\infty)$.  
\end{defin}

Note that the assumptions on $K$, $(\mu_t)_{t\geq 0}$ and $\phi$ ensure that
all the terms in (\ref{smolu}) and (\ref{flory}) make sense.

\medskip

Applying (\ref{smolu}) (or (\ref{flory})) with $\phi(x)=x$ 
(which does not belong to $C_c([0,\infty))$) would
clearly give $\Delta \phi =0$. Hence the total mass $\lc \mu_t(dx),x\rc$
is {\it a priori} constant as time evolves. 
However, 
for coagulation kernels satisfying $(A_\alpha)$, the gelation
phenomenon (that is, the loss of mass in finite time, or, equivalently,
the appearance of particles with infinite mass) is known to
occur \cite{emp,jeon}, which we recall now, together with other
properties. 

\begin{prop}\label{gege}
Consider a coagulation kernel $K$ satisfying $(A_\alpha)$ for some
$\alpha \in (0,1]$ and $\mu_0\in \mes$ such that $\left< \mu_0(dx),1\rc<\infty$
and $\lc \mu_0(dx),x\rc=1$. For any solution $(\mu_t)_{t\geq 0}$ of the
Smoluchowski or Flory equation, 
the gelation time
\begin{eqnarray}\label{dftgel}
T_{gel}:=\inf\{t\geq 0 \;\; : \;\;
\left<\mu_t(dx),x\right><\left<\mu_0(dx),x\right>\} 
\end{eqnarray}
is finite with the following upper estimate (here $c$ is defined in 
$(A_\alpha)$)
$$
T_{gel} \le \frac{\lc\mu_0(dx),x^{1-\alpha}\rc} {(1-2^{-\alpha})c}\,. 
$$

If $(\mu_t)_{t\geq 0}$ solves the Flory equation, then $t \mapsto \lc
\mu_t(dx), x \rc $ is continuous and strictly decreasing on
$(T_{gel},\infty)$, 
$$
\lim_{t \to \infty}\lc \mu_t(dx), x \rc=0 \;\;\mbox{ and }\;\;
\int_{T_{gel}+\e}^\infty \lc \mu_s(dx) , x^{1+\alpha} \rc ds <\infty 
$$
for all $\e>0$.
\end{prop}

The proof that gelation occurs is easier under $(A_\alpha)$ than the
general proof of Escobedo-Mischler-Perthame \cite{emp}, and we will
sketch it in the next section. This result expresses that particles
with infinite mass appear in finite time. Observe next that equations
$(S)$ and $(F)$ do not differ until gelation. The additional term in
equation $(F)$ represents the loss of finite particles with mass $x$,
proportionally to $l(x)$ and to the mass of the giant particles
$\left<\mu_0(dx)-\mu_s(dx), x \right>$.\\ 
Note that we are not able, and this is a well-known open problem, to
show that $t \mapsto \lc \mu_t(dx), x \rc $ is continuous at $t=T_{gel}$.

We finally consider a converging sequence of initial data.

\medskip

\textbf{Assumption $(I)$:} For each $n\in\nn\setminus\{0\}$, we are
given $\mu_0^n=m_n^{-1} \sum_{i=1}^n \delta_{x_i^n}$ for some $\left(
x_1^n,...,x_n^n \right)\in(0,\infty)^n$ and $m_n=x_1^n+...+x_n^n$. We
assume that there exists $\mu_0 \in \mes$ such that $\lc \mu_0(dx), x
\rc=1$ and $\lim_n \lc \mu^n_0(dx), \phi \rc =\lc \mu_0(dx), \phi \rc$
for all $\phi \in C_b([0,\infty))$, $C_b([0,\infty))$ denoting the
space of continuous and bounded functions on $[0,\infty)$. In
addition,  
$$
\lim_{\e \to 0}\sup_n \left< \mu_0^n(dx), \indiq_{(0,\e]} \right>=0\,.
$$

\medskip

We will actually not use explicitly all the assumptions in $(I)$ and
$(A_\alpha)$: some are just needed to apply the results of
\cite{fg}. We now state a compactness result which follows from
\cite{fg}. 

\begin{prop}\label{tension}
Consider a coagulation kernel $K$ satisfying $(A_\alpha)$ for some
$\alpha \in (0,1]$ and a sequence of initial conditions
$(\mu_0^n)_{n\ge 1}$ satisfying $(I)$. For each $a>0$ and $n\ge 1$, we
put $K_a:=K\ \indiq_{(0,a]\times (0,a]}$ and denote by
$(\mu^{n,a}_t)_{t\geq 0}$ the Marcus-Lushnikov process associated with
the pair $(K_a,\mu_0^n)$. The family $\{(\mu^{n,a}_t)_{t\geq
0}\}_{n\geq 1, a>0}$ is tight in $\mathbb{D}([0,\infty),\mes)$,
endowed with the Skorokhod topology associated with the vague topology
on $\mes$. 
\end{prop}

This proposition is proved in \cite[Theorem~2.3-i]{fg} (with the
choice of the subadditive function $\phi(x)=\sqrt{2C}(1+x)$, for which
$K(x,y)\leq \phi(x)\phi(y)$). Actually, it is stated in \cite{fg}
without the dependence on $a$, but the extension is straightforward.  

Notice here that, if $a\geq m_n$, the Marcus-Lushnikov process
$(\mu^{n,a}_t)_{t\geq 0}$ reduces to the standard Marcus-Lushikov
process associated with $(K,\mu_0^n)$. 

\medskip

We may finally state our main results. Recall that we assume
the total mass of the system to be initially $\lc \mu_0(dx),x\rc=1$.

\begin{theo}\label{result}
Consider a coagulation kernel $K$ satisfying $(A_\alpha)$ for some
$\alpha \in (0,1]$ and a sequence of initial conditions
$(\mu_0^n)_{n\ge 1}$ satisfying $(I)$. Consider also a sequence
$(a_n)_{n\geq 1}$ of positive real numbers 
such that $\lim_n a_n=\infty$. For each $n\geq 1$, let
$(\mu^{n,a_n}_t)_{t\geq 0}$ be the Marcus-Lushnikov process
associated with the pair $(K_{a_n},\mu_0^n)$ where $K_{a_n}:= K\
\indiq_{(0,a_n]\times (0,a_n]}$ and consider the weak limit
$(\mu_t)_{t\geq 0}$ in $\mathbb{D}([0,\infty),\mes)$ of a subsequence
$\left\{ (\mu^{n_k,a_{n_k}}_t)_{t\geq 0} \right\}_{k\geq 1}$. Then
$(\mu_t)_{t\geq 0}$ belongs a.s. to  $C([0,\infty),\mes)$ and enjoys
the following properties: 
\begin{enumerate}
\item Assume that $a_n=m_n$.\\ 
$(i)$ Then $(\mu_t)_{t\geq 0}$ solves a.s. the Flory equation with
coagulation kernel $K$ and initial condition $\mu_0$.\\ 
$(ii)$ Denote by $M_1^{n}(t)\geq M_2^{n}(t)\geq ... $ the ordered
sizes of the particles in the Marcus-Lushnikov
process $(\mu^{n,m_n}_t)_{t\geq 0}$, and define the ({\it a priori}
random) gelation time $T_{gel}$  of  $(\mu_t)_{t\geq 0}$ as in
(\ref{dftgel}). Then for all $\eta>0$ and $\beta>0$, 
\begin{eqnarray}
\lim_{k \to \infty} E \left[ \sup_{t\in[T_{gel}+\eta, \infty)} \left|
\frac{M_1^{n_k}(t)}{m_{n_k}} -\left(1-\left<\mu_t(dx), x\right>\right)
\right| \right] =0,\label{1i1} \\ 
\lim_{b \to \infty} \limsup_{k\to\infty} P
\left[\left(\int_{T_{gel}}^\infty \frac{1}{m_\nk}\sum_{i\geq 2}
M^{\nk}_i(s) \indiq_{[b,\infty)}\left(M^{\nk}_i(s) \right) ds
\right)\geq\beta \right] =0.\label{1i3} 
\end{eqnarray}
Furthermore, there is a positive constant $L$ depending only on $K$
such that, for all $\eta>0$ and $b>1$,  
\begin{eqnarray}
&\ds \lim_{k \to \infty} E \left[ \sup_{t\in[0,T_{gel}-\eta]}
\frac{M_1^{n_k}(t)}{m_{n_k}}\right] =0,\label{1i2} \\ 
&\ds \limsup_{k \to \infty} E \left[\int_0^{T_{gel}} \lc
\mu^{\nk,m_\nk}_s(dx),  x\indiq_{[b,\infty)}(x) \rc^2 ds \right]\leq
\frac{L}{b^\alpha}. \label{1i4} 
\end{eqnarray}
\item If $a_n/m_n\to 0$ as $n\to\infty$, then $(\mu_t)_{t\geq 0}$
solves a.s. the Smoluchowski equation with coagulation kernel $K$ and
initial condition $\mu_0$. 
\item If $a_n/m_n \to\gamma \in (0,1)$ as $n\to\infty$, then
$(\mu_t)_{t\in [0,T_1)}$ solves a.s. the Flory equation with
coagulation kernel $K$ and initial condition $\mu_0$ where 
\begin{equation}
T_{1}:=\inf\{ t > 0 \;\; : \;\; 1-\left<\mu_t(dx),x\right> \geq \gamma
\}. 
\end{equation}
\end{enumerate}
\end{theo}

Point~1-(i) is proved in \cite[Theorem~2.3-ii]{fg}. Remark that
(\ref{1i2}) is almost obvious while (\ref{1i4}) gives an estimate on
the tail of the mass distribution before gelation. The most
interesting estimate is of course (\ref{1i1}) which shows that, for $t>
T_{gel}$, the largest particle in the Marcus-Lushnikov process without
cut-off occupates a positive fraction of the total mass of the system
with a precise asymptotic. Finally, (\ref{1i3}) shows that, in some
sense, there is only one giant particle after gelation: the other
particles are rather small. Other results about the largest particles
for the kernel  
$$
K(x,y)=\frac{2(xy)^{1+\alpha}} {(x+y)^{1+\alpha}-x^{1+\alpha}-y^{1+\alpha}},
$$
which satisfies $(A_\alpha)$, were obtained by Aldous
\cite{aldousgp}. He however did not show that, after gelation, the
size of the largest particle is of order $\e m_n$. 

Point~2 seems to be new, and quite interesting. Indeed, we allow
arbitrary cut-off sequences $(a_n)$ which increase more slowly than
$(m_n)$. 

Finally, Point~3 can be explained in the following way: assume that
$a_n=\gamma m_n$ for all $n\ge 1$ and some $\gamma\in (0,1)$ and that
there is only one giant particle in $(\mu^{n,m_n}_t)_{t\geq 0}$. In
that situation, we then clearly have $(\mu^{n,\gamma m_n}_t)_{t\in
[0,T^n_1]}=(\mu^{n,m_n}_t)_{t \in [0, T_1^n]}$, where $T_1^n$ is the
first time at which the giant particle has a size greater than $\gamma
m_n$, i.e., it occupates a fraction $\gamma$ of the total mass of the
system.  Thus, $(\mu^{n,\gamma m_n}_t)_{t\in [0,T^n_1]}$ should
converge to $(\mu_t)_{t\in [0,T_1]}$, where $\mu$ solves the Flory
equation, and $T_1$ is the first time for which the giant particle
occupates a fraction $\gamma$ of the total mass in the Flory model. 

The proof of Theorem~\ref{result} is given in Section~\ref{theproof},
after establishing some properties of solutions to the Smoluchowski
and Flory coagulation equations in the next section. The final section
of the paper is devoted to numerical illustrations. 


\section{Properties of solutions to $(S)$ and $(F)$}\label{properties}
\setcounter{equation}{0}

Throughout this section, $K$ is a coagulation kernel satisfying
$(A_\alpha)$ for some $\alpha \in (0,1]$ and $\mu_0$ belongs to $\mes$
with total mass $\lc \mu_0(dx), x \rc=1$.  

\medskip

\begin{proof} {\bf of Proposition~\ref{gege}.}
Let $(\mu_t)_{t\geq 0}$ be a solution to the Smoluchowski equation
$(S)$ or the Flory equation $(F)$, and define $T_{gel}\in (0,\infty]$
by (\ref{dftgel}). Classical approximation arguments allow us to
use (\ref{smolu}) and (\ref{flory}) with
$\phi(x)=x^{1-\alpha}$. Indeed, it suffices to approximate $\phi$ by a
sequence of functions in $C_c([0,\infty))$ and to pass to the limit,
using the first inequality in
\begin{equation}
\label{spip}
\min(x,y)^{1-\alpha} \ge x^{1-\alpha} + y^{1-\alpha} -(x+y)^{1-\alpha}
\geq (2-2^{1-\alpha})\ \min(x,y)^{1-\alpha}\,, 
\end{equation}
which warrants that $K(x,y)|\Delta\phi(x,y)|\leq 2 C xy$ by $(A_\alpha)$.
We deduce from (\ref{smolu}), (\ref{flory}), and the second inequality
in (\ref{spip}) that, for all $t\geq 0$, 
\begin{eqnarray*}
\lc \mu_t(dx),x^{1-\alpha}\rc & \leq & \lc \mu_0(dx),x^{1-\alpha}\rc \\
&&- \frac{2-2^{1-\alpha}}{2}\intot \lc \mu_s(dx)\mu_s(dy), K(x,y) 
\min(x,y)^{1-\alpha} \rc \ ds.
\end{eqnarray*}
By virtue of (\ref{assumption}), $K(x,y) \min(x,y)^{1-\alpha} \geq c x
y$, whence 
\begin{equation}
\label{spirou}
(1-2^{-\alpha})c \intot \lc \mu_s(dx), x \rc^2 \ ds \leq \lc
\mu_0(dx),x^{1-\alpha}\rc 
\end{equation}
for all $t\geq 0$. Since $\lc \mu_s(dx), x \rc=\lc \mu_0(dx), x \rc=1$
for all $s\in [0,T_{gel})$, we realize that $T_{gel}$ has to be finite
for (\ref{spirou}) to hold true. A further consequence of
(\ref{spirou}) is that   
$$
(1-2^{-\alpha})c T_{gel} = (1-2^{-\alpha})c 
\int_0^{T_{gel}} \lc \mu_s(dx), x \rc^2\ ds 
\leq \lc\mu_0(dx),x^{1-\alpha}\rc, 
$$
whence
$$
T_{gel} \leq \frac{\lc\mu_0(dx),x^{1-\alpha}\rc} {(1-2^{-\alpha})c}\,.
$$
It also follows from (\ref{spirou}) that $t\longmapsto \lc \mu_t(dx),
x \rc$ belongs to $L^2(0,\infty)$ which, together with the
monotonicity and non-negativity of $t\longmapsto \lc \mu_t(dx), x \rc$
implies that $\lc \mu_t(dx), x \rc\longrightarrow 0$ as $t\to\infty$.  

\medskip

We now assume that  $(\mu_t)_{t\geq 0}$ solves the Flory equation
$(F)$ and prove that, for all $\e>0$,  
\begin{equation}\label{memfin}
\int_{T_{gel}+\e}^\infty \lc \mu_t(dx),x^{1+\alpha} \rc \ dt <\infty.
\end{equation}
To do so, we apply (\ref{flory}) with the choice $\phi_A(x)=\min
(x,A)$ for some positive real number $A$. Since $\Delta \phi_A$ is
non-positive, we get 
$$
\lc \mu_t,\phi_A\rc\leq \lc \mu_0,\phi_A\rc - \intot    
\lc \mu_s(dx),\min(x,A) l(x)\rc \left( 1 - \lc \mu_s(dx) ,x \rc
\right)\ ds. 
$$
Since $1=\lc \mu_0(dx) , x \rc \ge \lc \mu_s(dx), x \rc$ we may let
$A\to\infty$ and $t\to\infty$ in the above inequality and use the
Fatou lemma to deduce that  
\begin{equation}
\label{fantasio}
\int_0^\infty \lc \mu_s(dx) ,  xl(x) \rc \left( 1 - \lc \mu_s(dx) , x
\rc \right)\ ds \le 1. 
\end{equation}
Let $\e>0$. On the one hand, putting
$$
\delta_\e:= \inf_{t\geq T_{gel}+\e}\left\{ 1 - \lc \mu_t(dx),x \rc
\right\}, 
$$
it follows from the definition (\ref{dftgel}) of $T_{gel}$ that
$\delta_\e>0$. On the other hand, $x l(x)\ge c x^{1+\alpha}$ by
$(A_\alpha)$.  We therefore infer from (\ref{fantasio}) that  
$$
\int_{T_{gel}+\e}^\infty  \lc \mu_s(dx),x^{1+\alpha} \rc ds \leq 
\frac{1}{c\delta_\e}\,,
$$
whence (\ref{memfin}).

We now check that $t\mapsto \lc \mu_t(dx),x\rc$ is continuous on $(T_{gel},
\infty)$.
Using once more (\ref{flory}) with the choice $\phi_A(x)=\min(x,A)$,
we obtain for $T_{gel}<s<t$ 
\begin{eqnarray*}
\left< \mu_t-\mu_s, \phi_A \right> &=& \frac{1}{2} \int_s^t
\left< \mu_\tau(dx) \mu_\tau (dy), K(x,y)\Delta\phi_A(x,y) \right>
d\tau \\ 
&&- \int_s^t  \left< \mu_\tau(dx), \phi_A(x) l(x) \right> 
\left( 1 - \lc \mu_\tau(dx), x \rc \right) d\tau. 
\end{eqnarray*}
Clearly $\Delta\phi_A(x,y)\to 0$ as $A\to\infty$ for all $(x,y)\in
(0,\infty)^2$ while $(A_\alpha)$ warrants that  
$$
K(x,y)|\Delta\phi_A(x,y)|\leq C(x^\alpha y + x y^\alpha) \min(x,y)
\leq C (x^{1+\alpha}y + x y^{1+\alpha}). 
$$
Using (\ref{memfin}) and the Lebesgue dominated convergence theorem,
we obtain 
\begin{eqnarray}\label{exmass}
\lc \mu_t(dx)-\mu_s(dx),x \rc = - \int_s^t \lc\mu_\tau(dx),x l(x)\rc
\left( 1 - \lc \mu_\tau(dx),x \rc \right) d\tau\,.
\end{eqnarray}
Using again (\ref{memfin}) and that $| \lc\mu_\tau(dx),x l(x)\rc
\left( 1 - \lc \mu_\tau(dx),x \rc \right)| \leq C  
\lc\mu_\tau(dx),x^{1+\alpha}\rc$ by $(A_\alpha)$, 
we conclude that $t \longmapsto \lc \mu_t(dx),x \rc$ is
continuous on $(T_{gel},\infty)$.  

It remains to check that $t \longmapsto \lc \mu_t(dx),x \rc$ is
strictly decreasing on $t\in (T_{gel},\infty)$. According to
(\ref{exmass}) this is true as long as $\mu_\tau\ne 0$ for
$T_{gel}<s<\tau<t$: it thus suffices to show that $\mu_t\ne 0$ for
all $t\ge 0$. For that purpose, we take $\phi(x)=x\indiq_{(0,A]}(x)$
in (\ref{flory}) where $A>0$ is chosen so that $\lc \mu_0(dx),
x\indiq_{(0,A]}(x)\rc>0$ (such an $A$ always exists as
$\lc\mu_0(dx),x\rc=1$). Thanks to (\ref{assumption}), we get 
\begin{eqnarray}
&& \frac{d}{dt} \lc \mu_t,\phi\rc \geq - \lc\mu_t(dx)\mu_t(dy),
K(x,y) x\indiq_{(0,A]}(x) \rc - \lc \mu_t(dx), x l(x)
\indiq_{(0,A]}(x) \rc \nonumber\\ 
&&\hskip 1cm \geq - C  \lc \mu_t(dx)\mu_t(dy),
(x^{\alpha+1}y+x^2y^\alpha) \indiq_{(0,A]}(x)\rc -C \lc \mu_t(dx),
x^{1+\alpha} \indiq_{(0,A]}(x)\rc\nonumber\\ 
&&\hskip 1cm \geq  -C (A^\alpha+A) \lc \mu_t, \phi \rc \lc
\mu_t(dx),x+x^\alpha \rc - C A^\alpha  \lc \mu_t, \phi \rc.
\end{eqnarray}
Since $t\longmapsto \lc\mu_t(dx), x\rc$ and $t\longmapsto
\lc\mu_t(dx), x^\alpha\rc$ are non-increasing and $\lc
\mu_0(dx),x+x^\alpha \rc < \infty$, we conclude that  
$$
\frac{d}{dt} \lc \mu_t,\phi\rc \geq -C_A \left( 1 + \lc \mu_0(dx) ,
x^\alpha \rc \right)\ \lc \mu_t,\phi\rc  \geq -C_{A,\mu_0} \lc \mu_t,\phi\rc
$$
for all $t\ge 0$ for some constant $C_{A,\mu_0}>0$. 
Consequently, $\lc \mu_t, \phi \rc>0$ for all
$t\ge 0$ as the choice of $A$ warrants that $\lc \mu_0, \phi
\rc>0$, and the proof of Proposition~\ref{gege} is complete. 
\end{proof}

\medskip

Next, as a preliminary step towards the proof of Theorem~\ref{result}
Point~1-(ii), we show that solutions to the Smoluchowski and Flory
coagulation equations do not coincide after the gelation time.  

\begin{cor}\label{flnesm}
Let $(\mu_t)_{t\geq 0}$ and $(\nu_t)_{t\geq 0}$ be solutions to the
Smoluchowski equation $(S)$ and the Flory equation $(F)$,
respectively, (with the same coagulation kernel $K$ and initial
condition $\mu_0$), and assume further their respective gelation times
coincide, that is, 
$$
T_{gel}:=\inf\{t\geq 0 \;\; : \;\;
\left<\mu_t(dx),x\right><\left<\mu_0(dx),x\right>\} = \inf\{t\geq 0
\;\; : \;\; \left<\nu_t(dx),x\right><\left<\mu_0(dx),x\right>\}\,. 
$$
Then, for each $\e>0$, there exists $s_\e\in(T_{gel},T_{gel}+\e)$ such
that $\mu_{s_\e} \ne \nu_{s_\e}$. 
\end{cor}

\begin{proof}\textbf{.}
Consider $\e>0$.\\
Either $t\longmapsto \lc \mu_t(dx),x^{1+\alpha} \rc$ does not belong
to $L^1\left( T_{gel}+(\e/2), T_{gel}+\e \right)$ and $\mu_t$ cannot
coincide with $\nu_t$ on $\left( T_{gel}+(\e/2), T_{gel}+\e \right)$
since $t\longmapsto \lc \nu_t(dx),x^{1+\alpha} \rc$ 
belongs to $L^1\left( T_{gel}+(\e/2), T_{gel}+\e
\right)$ by Proposition~\ref{gege}.\\ 
Or $t\longmapsto \lc \mu_t(dx),x^{1+\alpha} \rc$ belongs to $L^1\left(
T_{gel}+(\e/2), T_{gel}+\e \right)$ and it is not difficult to check
that this property and (\ref{smolu}) entail that
$\left<\mu_t(dx),x\right> = \left<\mu_{T_{gel}+(\e/2)}(dx),x\right>$ for
$t\in \left[ T_{gel}+(\e/2), T_{gel}+\e \right]$: indeed, take 
$\phi_A(x)=\min (x, A)$ in (\ref{smolu}) and pass to the limit as 
$A\to \infty$ using that $\Delta \phi_A(x,y) \to 0$ and the time 
integrability of $t\longmapsto \lc \mu_t(dx),x^{1+\alpha} \rc$.
Owing to the strict
monotonicity of $t\longmapsto \lc \nu_t(dx),x \rc$ established in
Proposition~\ref{gege}, the previous property of $\mu_t$ excludes that
$\mu_t=\nu_t$ for all $t\in \left( T_{gel}+(\e/2), T_{gel}+\e \right)$
and completes the proof of Corollary~\ref{flnesm}. 
\end{proof}


\section{Proof of the main results}\label{theproof}\setcounter{equation}{0}

We fix a coagulation kernel $K$ satisfying $(A_\alpha)$ for some
$\alpha \in (0,1]$ and a sequence of initial data $(\mu_0^n)_{n\ge 1}$
satisfying $(I)$. Next, for $a>0$ and $n\ge 1$, we put
$K_a(x,y):=K(x,y)\indiq_{(0,a]\times(0,a]}$ and denote by
$(\mu^{n,a}_t)_{t\geq 0}$ the Marcus-Lushnikov process associated with
the pair $(K_a,\mu_0^n)$. According to Definition~\ref{dfml} we may
write  
\begin{equation}
\label{gaston}
\mu^{n,a}_t = \frac{1}{m_n}\sum_i \delta_{M_i^{n,a}(t)} \;\;\mbox{
with }\;\; M_1^{n,a}(t)\geq M_2^{n,a}(t)\geq M_3^{n,a}(t) \geq ... 
\end{equation}
for all $t\ge 0$, $n\ge 1$, and $a>0$.

We next recall that the space $\mathbb{D}([0,\infty),\mes)$ is endowed
with the Skorokhod topology associated with the vague convergence
topology on $\mes$ (see Ethier-Kurtz \cite{ek} for further
information), and denote by $d$ a distance on $\mes$ metrizing the
vague convergence topology.  

Marcus-Lushnikov processes have some martingale properties, which are
immediately obtained from (\ref{gml}), see also
\cite[Section~4]{norris}. 

\begin{lemma}\label{martingale}
For all $\phi\in L^\infty_{loc}(0,\infty)$ and $t\geq 0$, we have
\begin{eqnarray}
&&\lc \mu^{n,a}_t, \phi \rc = \lc \mu^{n}_0, \phi \rc +
O^{n,a}_t(\phi) \nonumber\\
&&\hskip0.5cm + \frac{1}{2 m_n^2}\intot \sum_{i\ne j}
K_a(M_i^{n,a}(s),M_j^{n,a}(s)) \Delta\phi(M_i^{n,a}(s),M^{n,a}_j(s))\
ds \nonumber\\
&&= \lc \mu^{n}_0, \phi \rc + O^{n,a}_t(\phi)+ \frac{1}{2}\intot
\lc  \mu^{n,a}_s(dx)\mu^{n,a}_s(dy), K_a(x,y)\Delta \phi(x,y) \rc\ ds
\nonumber\\ 
&&\hskip0.5cm - \frac{1}{2m_n} \intot \lc \mu^{n,a}_s(dx),
K_a(x,x)\Delta\phi(x,x)\rc\ ds \label{eqml}
\end{eqnarray}
where $\Delta \phi$ is defined in (\ref{dfdelta}), and $O^{n,a}(\phi)$
is a martingale starting from $0$ with (predictable) quadratic
variation  
\begin{eqnarray*}
\lc O^{n,a}(\phi) \rc_t &=& \frac{1}{2m_n}\intot \lc
\mu^{n,a}_s(dx)\mu^{n,a}_s(dy), K_a(x,y) \left[
\Delta\phi(x,y)\right]^2 \rc\ ds \\ 
&&-\frac{1}{2 m_n^2}\intot \lc  \mu^{n,a}_s(dx), K_a(x,x) \left[
\Delta\phi(x,x)\right]^2 \rc\ ds. 
\end{eqnarray*}
Furthermore, if $\phi:(0,\infty)\to\rr$ is a subadditive function,
that is, $\phi(x+y)\leq \phi(x)+\phi(y)$ for $(x,y)\in (0,\infty)^2$,
then $t\mapsto \left<\mu^{n,a}_t,\phi \right>$ is a.s. a
non-increasing function. 
\end{lemma}

We carry on with some easy facts.

\begin{lemma}\label{toujoursok}
Let $(a_n)_{n\geq 1}$ be a sequence of positive numbers. Then any weak
limit $(\mu_t)_{t\geq 0}$ of the sequence $\left\{
(\mu^{n,a_n}_t)_{t\geq 0} \right\}_{n\ge 1}$ belongs a.s. to
$C([0,\infty),\mes)$, and both $t\mapsto \lc \mu_t(dx),x \rc$ and
$t\mapsto \lc \mu_t(dx),1 \rc$ are a.s. non-increasing functions. 
Furthermore, 
\begin{eqnarray}\label{kappa}
\sup_{n\ge 1} \sup_{t\geq 0} \lc \mu^{n,a_n}_t(dx), 1+x \rc = \kappa
:= \sup_n \lc \mu^{n}_0(dx), 1+x \rc < \infty, 
\end{eqnarray}
and for all $\phi\in C_c([0,\infty))$ and $T>0$,
\begin{eqnarray}
&&\lim_{n\to\infty} E\left[ \sup_{t\in [0,T]} \left| \frac{1}{2m_n}
\intot \lc \mu^{n,a_n}_t(dx), K_{a_n}(x,x) \Delta\phi(x,x) \rc\ ds
\right|\right] =0, 
\label{diagz} \\
&&\lim_{n\to\infty} E\left[ \sup_{t\in [0,T]}
\left(O^{n,a_n}_t(\phi)\right)^2\right] = 0. 
\label{martoz}
\end{eqnarray}
\end{lemma}

\begin{proof}\textbf{.} First, if $\phi \in C_b([0,\infty))$, the
jumps of $\lc \mu^{n,a_n}_t,\phi\rc$ are of the form $m_n^{-1}
\Delta\phi(x,y)$ and clearly converge to zero as $n\to\infty$ since
$m_n\to\infty$. Hence any weak limit $(\mu_t)_{t\geq 0}$ belongs to
$C([0,\infty),\mes)$ a.s.  

Consider next a family $(\xhi_b)_{b>0}$ of continuous non-increasing
functions on $(0,\infty)$ such that $\xhi_b(x)=1$ for $x\leq b$ and
$\xhi_b(x)=0$ for $x\geq b+1$, and a non-negative subadditive function
$\phi$. Then, on the one hand, $\phi \xhi_b$ is also subadditive and
Lemma~\ref{martingale} ensures that $t\longmapsto\lc
\mu^{n,a_n}_t, \phi\xhi_b\rc$ is a.s. non-increasing for all
$n\geq 1$. On the other hand, since $\phi \xhi_b\in C_c([0,\infty))$,
it follows from the definition of $(\mu_t)_{t\ge 0}$ that there is a
subsequence $(\nk)_{k\ge 1}$, $\nk\to\infty$, such that $\left\{
\left( \lc \mu^{\nk,a_\nk}_t, \phi\xhi_b\rc \right)_{t\geq 0}
\right\}_{k\ge 1}$ converges in law towards $(\lc \mu_t,
\phi\xhi_b\rc)_{t\geq 0}$ for each fixed $b>0$ as
$k\to\infty$. Therefore, $t\longmapsto \lc \mu_t, \phi\xhi_b\rc$
is a.s. non-increasing for each $b>0$. Since $\left( \lc \mu_t,
\phi\xhi_b\rc \right)_{b>0}$ converges to $\lc \mu_t, \phi \rc$ as
$b\to\infty$ for each $t\ge 0$, we conclude that $t\mapsto \lc
\mu_t, \phi \rc$ is a.s. non-increasing. Applying this result to
$\phi(x)=1$ and $\phi(x)=x$, we obtain that both $t\mapsto \lc
\mu_t(dx),x \rc$ and $t\mapsto \lc \mu_t(dx),1 \rc$ are
a.s. non-increasing functions of time.  

Next, since $x\mapsto 1+x$ is subadditive, Lemma~\ref{martingale}
implies that we have a.s. $\lc \mu^{n,a_n}_t(dx), 1+x \rc \le \lc
\mu^{n}_0(dx), 1+x \rc$ for $n\geq 1$ and $t\ge 0$, and $\lc
\mu^{n}_0(dx), 1+x \rc$ is bounded uniformly with respect to $n$ by
assumption $(I)$. 

Consider finally $\phi \in C_c([0,\infty))$ with support included in
$[0,R]$ for some $R>0$. By (\ref{assumption}),
$|K_{a_n}(x,x)\Delta\phi(x,x)|\leq 6C \|\phi\|_{L^\infty}\
R^{1+\alpha}$, whence 
$$
\left|\lc \mu^{n,a_n}_t(dx), K_{a_n}(x,x)\Delta\phi(x,x) \rc\right|
\leq 6C\kappa\|\phi\|_{L^\infty}\ R^{1+\alpha} \;\;\mbox{ a.s. } 
$$
by (\ref{kappa}), from which (\ref{diagz}) readily follows since
$m_n\to\infty$. By a similar argument, we establish that $E \left[\lc
O^{n,a_n}(\phi) \rc_t \right]\longrightarrow 0$ as $n\to\infty$, which
implies (\ref{martoz}) by Doob's inequality. 
\end{proof}

\medskip

We now prove a fundamental estimate which provides a control on the
large masses contained in $\mu^{n,a}_t$. 

\begin{lemma}\label{fund}
There exists a positive real number $L$ depending only on $c$ and
$\alpha$ in $(A_\alpha)$ such that 
\begin{eqnarray}\label{tiptop}
E\left[\int_0^\infty \frac{1}{m_n^2} \sum_{i\ne j }
M_i^{n,a}(s)M_j^{n,a}(s) \indiq_{[b,a]}\left( M_i^{n,a}(s) \right)
\indiq_{[b,a]}\left( M_j^{n,a}(s) \right)\ ds \right] \leq
\frac{L}{b^\alpha}
\end{eqnarray}
for all $n\geq 1$, $a>0$, and $b\in (0,a)$, the $M_i^{n,a}$ being
defined in (\ref{gaston}).  
\end{lemma}

\begin{proof}\textbf{.}
To prove this estimate, we use (\ref{eqml}) with $\phi(x)=x^{1-\alpha}
\min{(x,b)}^\alpha$ for some $b\in (0,a)$. We first notice that $\lc
\mu_0^n, \phi \rc \leq \lc \mu^n_0(dx),x\rc = 1$ and $\lc
\mu^{n,a}_t, \phi \rc \geq 0$ for all $t\geq 0$ and $n\ge 1$. In
addition, $\phi$ is subadditive so that $\Delta \phi(x,y)$ is always
non-positive and we infer from (\ref{assumption}) and (\ref{spip})
that  
$$
K_a(x,y)\Delta \phi (x,y) \leq -  (2-2^{1-\alpha})c b^\alpha  xy
\indiq_{[b,a]}(x) \indiq_{[b,a]}(y) 
$$
for $(x,y)\in (0,\infty)^2$. Taking expectations in (\ref{eqml}) and
using the above inequalities, we obtain  
$$
0  \leq 1 - \frac{b^{\alpha}}{L} E\left[ \intot \frac{1}{m_n^2} 
\sum_{i\ne j } M_i^{n,a}(s)M_j^{n,a}(s) \indiq_{[b,a]}\left(
M_i^{n,a}(s) \right) \indiq_{[b,a]}\left( M_j^{n,a}(s) \right)\ ds
\right] 
$$
for all $t\geq 0$, with $1/L:=c(1-2^{-\alpha})$. We conclude the
proof by letting $t\to\infty$ in the previous inequality. 
\end{proof}

\medskip

We now turn to the proof of Theorem~\ref{result} and first recall that
Point~1-(i) is included in \cite[Theorem~2.3-ii]{fg} as
$(\mu^{n,m_n}_t)_{t\geq 0}$ is the standard Marcus-Lushnikov process
associated with $(K,\mu_0^n)$.  

\medskip

\begin{proof} {\bf of Point~2 of Theorem~\ref{result}.} 
Let $(a_n)_{n\geq 1}$ be a sequence of positive real numbers
satisfying $a_n\to\infty$ and $a_n/m_n\to 0$ as $n\to\infty$. We
consider the limit $(\mu_t)_{t\geq 0}$ of a subsequence $\left\{
(\mu^{n_k,a_{n_k}}_t)_{t\geq 0} \right\}_{k\ge 1}$ in the sense that
a.s. 
\begin{equation}
\label{prunelle}
\lim_{k\to\infty} \sup_{[0,T]} d(\mu^{n_k,a_{n_k}}_t,\mu_t)=0
\;\;\mbox{ for all }\;\; T>0\,, 
\end{equation}
the existence of such a limit being guaranteed by
Proposition~\ref{tension}, Lemma \ref{toujoursok} 
(which ensures the time continuity of $\mu_t$) 
and the Skorokhod representation theorem. We
now aim at showing that $(\mu_t)_{t\geq 0}$ solves a.s. the
Smoluchowski equation $(S)$ and proceed in two steps. 

\noindent {\it Step 1.} We first deduce from Lemma~\ref{fund} that 
\begin{equation}\label{ccc}
E \left[\int_0^T \lc \mu_t^{n,a_n}(dx), x \indiq_{[b,a_n]}(x) \rc^2\
dt \right] \leq \left( \frac{L}{b^\alpha} + T \frac{a_n}{m_n} \right) 
\end{equation}
for all $b>0$, $n\geq 1$, and $T>0$. Indeed, we have a.s., for all
$t\geq 0$, 
\begin{eqnarray}
&&\frac{1}{m_n^2} 
\sum_{i\ne j } M_i^{n,a_n}(t)M_j^{n,a_n}(t)
\indiq_{[b,a_n]}\left(M_i^{n,a_n}(t)\right)
\indiq_{[b,a_n]}\left(M_j^{n,a_n}(t)\right) \nonumber\\ 
&&= \lc \mu^{n,a_n}_t(dx) \mu^{n,a_n}_t(dy), xy \indiq_{[b,a_n]}(x)
\indiq_{[b,a_n]}(y) \rc - \frac{1}{m_n} \lc \mu^{n,a_n}_t(dx), x^2
\indiq_{[b,a_n]}(x) \rc \nonumber\\
&&\geq \lc \mu^{n,a_n}_t(dx) , x \indiq_{[b,a_n]}(x) \rc^2 -
\frac{a_n}{m_n} \lc \mu^{n,a_n}_t(dx), x \rc
\label{ttoopp} \\ 
&& \ge \lc \mu^{n,a_n}_t(dx) , x \indiq_{[b,a_n]}(x) \rc^2 -
\frac{a_n}{m_n}\,, \nonumber 
\end{eqnarray}
hence (\ref{ccc}) after integrating over $(0,T)$, taking expectation,
and using Lemma~\ref{fund} (with $a=a_n$). 

\smallskip

{\it Step 2.} By Lemma~\ref{toujoursok}, we already know that
$t\mapsto \lc \mu_t(dx), x \rc$ and $t\mapsto \lc \mu_t(dx), 1 \rc$
are a.s. non-increasing functions. Consider now $\phi \in C_c
([0,\infty))$. The convergence (\ref{prunelle}) and the assumption
$(I)$ ensure that $\lc \mu^{\nk,a_\nk}_t, \phi \rc \longrightarrow \lc
\mu_t,\phi \rc$ a.s. for all $t\geq 0$ and $\lc \mu^{\nk}_0, \phi \rc
\longrightarrow \lc \mu_0,\phi \rc$ as $k\to\infty$. Recalling
(\ref{eqml}), (\ref{diagz}), and (\ref{martoz}), we realize that
$(\mu_t)_{t\geq 0}$ solves (\ref{smolu}) provided that we check 
that $B_k(t) \longrightarrow B(t)$ (for instance in $L^1$)
as $k\to\infty$ for all $t\geq 0$, where 
\begin{eqnarray*}
B_k(t)&:=& \intot \lc \mu^{\nk,a_\nk}_s(dx)\mu^{\nk,a_\nk}_s(dy), 
K_{a_\nk}(x,y)\Delta\phi(x,y) \rc\ ds, \\
B(t)  &:=& \intot \lc \mu_s(dx)\mu_s(dy), K(x,y)\Delta\phi(x,y) \rc\ ds.
\end{eqnarray*}
For that purpose, we consider a family $(\xhi_b)_{b>0}$ of continuous
non-increasing functions on $[0,\infty)$ such that $\xhi_b(x)=1$ for
$x\in (0,b]$ and  $\xhi_b(x)=0$ for $x\in[b+1,\infty)$, and put 
\begin{eqnarray*}
B_k(t,b)&:=& \intot \lc \mu^{\nk,a_\nk}_s(dx)\mu^{\nk,a_\nk}_s(dy),
K_{a_\nk}(x,y) \Delta\phi(x,y) \xhi_b(x) \xhi_b(y) \rc\ ds, \\ 
B(t,b)  &:=& \intot \lc \mu_s(dx)\mu_s(dy),  K(x,y) \Delta\phi(x,y)
\xhi_b(x) \xhi_b(y) \rc\ ds. 
\end{eqnarray*}
On the one hand, it follows from $(A_\alpha)$, the boundedness of
$\phi$, the bounds $x^\alpha \leq 1+x$ and 
$$
\sup_{t\geq 0} \lc \mu_t(dx),1+x \rc=\lc \mu_0(dx),1+x \rc\,,
$$
and the Lebesgue dominated convergence theorem that
\begin{equation}
\label{pim}
\lim_{b\to\infty} E[|B(t,b)-B(t)|]=0.
\end{equation}
On the other hand, for each $b\in (0,\infty)$, we have $K_{a_\nk}(x,y)
\xhi_b(x)\xhi_b(y)= K(x,y)\xhi_b(x)\xhi_b(y)$ for all $(x,y) \in
(0,\infty)^2$ as soon as $b+1 \leq a_\nk$, the latter being true for
$k$ sufficiently large. Consequently, since $(x,y)\longmapsto
K(x,y)\xhi_b(x)\xhi_b(y)\Delta\phi(x,y)$ belongs to
$C_c([0,\infty)^2)$, the convergence (\ref{prunelle}) entails that
$B_k(t,b) \longrightarrow B(t,b)$ for all $t\ge 0$ a.s. as
$k\to\infty$. Thanks to $(A_\alpha)$ and (\ref{kappa}) we may apply
the Lebesgue dominated convergence theorem to obtain that 
\begin{equation}
\label{pam}
\lim_{k\to\infty} E[|B_k(t,b)-B(t,b)|]=0 \;\;\mbox{ for each }\;\; b>0\,.
\end{equation}
Finally, owing to (\ref{assumption}), we have for $k$ sufficiently
large (such that $a_\nk \geq b$)  
\begin{eqnarray*}
&& E[|B_k(t,b)-B_k(t)|] \\
&& \le 3 \|\phi\|_{L^\infty} E\left[ \intot \lc
\mu^{n_k,a_\nk}_s(dx)\mu^{n_k,a_\nk}_s(dy), K_{a_\nk}(x,y) \left( 1
-\xhi_b(x) \xhi_b(y) \right) \rc\ ds  \right]\\ 
&& \leq 6C \|\phi\|_{L^\infty} E\left[ \intot \lc
\mu^{n_k,a_\nk}_s(dx)\mu^{n_k,a_\nk}_s(dy), x^\alpha y
\indiq_{(0,a_\nk]}(x) \indiq_{(0,a_\nk]}(y) \left( 1 -\xhi_b(x)
\xhi_b(y) \right) \rc\ ds  \right]\\ 
&& \leq 6C \|\phi\|_{L^\infty} E\left[ \intot \lc
\mu^{\nk,a_\nk}_s(dx)\mu^{\nk,a_\nk}_s(dy), x^\alpha y
\indiq_{(0,a_\nk]}(x) \indiq_{[b,a_\nk]}(y) \rc\ ds \right]\\ 
&& + 6C \|\phi\|_{L^\infty} E\left[ \intot \lc
\mu^{\nk,a_\nk}_s(dx)\mu^{\nk,a_\nk}_s(dy), x^\alpha y
\indiq_{[b,a_\nk]}(x) \indiq_{(0,a_\nk]}(y) \rc\ ds \right]\\ 
&&\leq 6C \|\phi\|_{L^\infty} E\left[ \intot
\lc \mu^{\nk,a_\nk}_s(dx), (1+x) \indiq_{(0,a_\nk]}(x) \rc \lc
\mu^{\nk,a_\nk}_s(dy), y \indiq_{[b,a_\nk]}(y)\rc\ ds \right]\\ 
&&+ \frac{6C \|\phi\|_{L^\infty}}{b^{1-\alpha}} E\left[ \intot \lc
\mu^{\nk,a_\nk}_s(dx), x \indiq_{[b,a_\nk]}(x) \rc \lc
\mu^{\nk,a_\nk}_s(dy) ,y \indiq_{(0,a_\nk]}(y) \rc\ ds \right]. 
\end{eqnarray*}
We then infer from (\ref{kappa}), the Cauchy-Schwarz inequality, and
(\ref{ccc}) that, for $b\geq 1$,
\begin{eqnarray*}
E[|B_k(t,b)-B_k(t)|]
&\leq&  12\kappa C \|\phi\|_{L^\infty}
E\left[ \intot \lc \mu^{\nk,a_\nk}_s(dx), x \indiq_{[b,a_\nk]}(x) \rc\
ds \right] \\ 
&\leq & 12\kappa C \|\phi\|_{L^\infty} t^{1/2} E\left[ \intot \lc
\mu^{\nk,a_\nk}_s(dx), x \indiq_{[b,a_\nk]}(x) \rc^2\ ds
\right]^{1/2}\\ 
&\leq & 12\kappa C \|\phi\|_{L^\infty} t^{1/2}
\left(\frac{L}{b^{\alpha}} + t\frac{a_\nk}{m_\nk}\right)^{1/2}.  
\end{eqnarray*}
Since $a_n/m_n\to 0$ as $n\to\infty$, we may first let $k\to\infty$
and then $b\to\infty$ in the previous inequality to conclude that  
\begin{equation}\label{poum}
\lim_{b \to \infty} \limsup_{k\to\infty} E[|B_k(t,b)-B_k(t)|] =0.
\end{equation}
Combining (\ref{pim}), (\ref{pam}), and (\ref{poum}) ends the proof.
\end{proof}

\medskip

We next complete the proof of Point~1 of Theorem~\ref{result}.

\medskip

\begin{proof} {\bf of Point~1-$(ii)$ of Theorem~\ref{result}.}
Recall that we are in the situation where $a_n=m_n$, so that
$(\mu_t^{n,a_n})_{t\ge 0}=(\mu_t^{n,m_n})_{t\ge 0}$ is the classical
Marcus-Lushnikov process associated with $(K,\mu_0^n)$ for each $n\ge
1$. Let $(\mu_t)_{t\geq 0}$ be the limit of a subsequence $\left\{
(\mu^{n_k,m_{n_k}}_t)_{t\geq 0} \right\}_{k\ge 1}$ in the sense that
a.s. 
\begin{equation}
\label{marsu}
\lim_{k\to\infty} \sup_{[0,T]} d(\mu^{\nk,m_\nk}_t,\mu_t) = 0
\;\;\mbox{ for all }\;\; T>0\,,
\end{equation}
the existence of such a limit following by the same arguments 
as (\ref{prunelle}).

We already know from \cite{fg} that $(\mu_t)_{t\geq 0}$ solves
a.s. the Flory equation. We define the ({\it a priori} random) gelling
time $T_{gel}$ of $(\mu_t)_{t\geq 0}$ by (\ref{dftgel}) and write  
\begin{equation}
\label{repf}
\mu^{n,m_n}_t = \frac{1}{m_n}\sum_i \delta_{M_i^n(t)} \;\;\mbox{ with
}\;\; M_1^n(t)\geq M_2^n(t)\geq M_3^n(t) \geq ... 
\end{equation}
for all $t\ge 0$ and $n\ge 1$.

\medskip

As before we denote by $(\xhi_b)_{b>0}$ a family of continuous
non-increasing functions such that $\xhi_b(x)=1$ for $x\in [0,b]$ and
$\xhi_b(x)=0$ for $x\geq b+1$. We start with the proof of (\ref{1i2})
which is almost immediate. Since $t \mapsto M_1^n(t)/m_n$ is
a.s. non-decreasing and bounded by $1$, it suffices to check that
$P[M^\nk_1(T_{gel}-\eta)\geq \delta m_\nk]\longrightarrow 0$ as
$k\to\infty$ for all $\eta>0$ and $\delta>0$. For that purpose, fix
$b>0$. Since $1-\xhi_b\ge 0$, it follows from (\ref{repf}) that  
\begin{eqnarray*}
\lc \mu^{n_k,m_\nk}_t(dx), x (1-\xhi_b(x)) \rc & \ge &
\frac{M_1^{n_k,m_\nk}(t) \left( 1 - \xhi_b\left( M_1^{n_k,m_\nk}(t)
\right) \right)}{m_\nk}, \\ 
1 - \lc \mu^{n_k,m_\nk}_t(dx), x \xhi_b(x) \rc & \ge &
\frac{M_1^{n_k,m_\nk}(t)}{m_\nk} \indiq_{[b+1,\infty)}\left(
M_1^{n_k,m_\nk}(t) \right)\,. 
\end{eqnarray*}
For $k$ large enough we have $\delta m_{n_k}\ge b+1$ and thus
\begin{eqnarray*}
E\left[ 1 - \lc \mu^{n_k,m_\nk}_{T_{gel}-\eta}(dx), x \xhi_b(x) \rc
\right] & \ge & E\left[ \frac{M_1^{n_k,m_\nk}(T_{gel}-\eta)}{m_\nk}
\indiq_{[\delta m_{n_k},\infty)}\left( M_1^{n_k,m_\nk}(T_{gel}-\eta)
\right) \right] \\ 
& \ge & \delta P\left[ M^\nk_1(T_{gel}-\eta)\geq \delta m_\nk
\right]\,. 
\end{eqnarray*}
Now, thanks to the compactness of the support of $\xhi_b$ and
(\ref{marsu}), the sequence 
$$
\left( \lc
\mu^{n_k,m_\nk}_{T_{gel}-\eta}(dx), x \xhi_b(x) \rc \right)_{k\ge 1}
$$
converges a.s. to $\lc \mu_{T_{gel}-\eta}(dx), x \xhi_b(x) \rc$ and is
bounded by (\ref{kappa}). We may then let $k\to\infty$ in the above
inequality to obtain 
$$
E\left[ 1 - \lc \mu_{T_{gel}-\eta}(dx), x \xhi_b(x) \rc \right] \ge
\delta \limsup_{k\to\infty} P\left[ M^\nk_1(T_{gel}-\eta)\geq \delta
m_\nk \right]\,. 
$$ 
Next, owing to the definition of $T_{gel}$, the sequence $\left( \lc
\mu_{T_{gel}-\eta}(dx), x \xhi_b(x) \rc \right)_{b>0}$ converges
towards $\lc \mu_{T_{gel}-\eta}(dx), x \rc =1$ as $b\to\infty$ and is
bounded by $1$. Passing to the limit as $b\to\infty$ in the previous
inequality entails that $P\left[ M^\nk_1(T_{gel}-\eta)\geq \delta
m_\nk \right]\longrightarrow 0$ as $k\to\infty$, which is the claimed
result. The limit (\ref{1i2}) then follows.  

\medskip

We now turn to the proof of (\ref{1i4}). Let $b>0$. Since
$a_n=m_n\to\infty$ as $n\to\infty$, we have $a_{n_k}>b$ for $k$ large
enough and it follows from Lemma~\ref{fund} (with $a=a_n)$ by an
argument similar to (\ref{ttoopp}) that
(recall that all the particles represented in $\mu^{n_k,m_\nk}_t$
are smaller than $m_\nk$ by construction, see Definition \ref{dfml})
\begin{eqnarray}
&&E\left[\int_0^{T_{gel}} \lc \mu^{n_k,m_\nk}_t(dx), x
\indiq_{[b,\infty)}(x) \rc^2\ dt \right] \nonumber\\ 
&&\leq \frac{L}{b^\alpha} + E\left[\int_0^{T_{gel}}\frac{1}{m_\nk}\lc
\mu^{n_k,m_\nk}_t(dx), x^2 \indiq_{[b,\infty)}(x) \rc\ dt \right]
\nonumber\\ 
&&\leq \frac{L}{b^\alpha} +
E\left[\int_0^{T_{gel}}\frac{M_1^{n_k}(t)}{m_\nk}\lc
\mu^{n_k,m_\nk}_t(dx), x \rc\ dt \right] \nonumber\\ 
&&\leq \frac{L}{b^\alpha} + E\left[\int_0^{T_{gel}}
\frac{M_1^{n_k}(t)}{m_\nk}\ dt \right]\,. \label{chch} 
\end{eqnarray}
Since $M^{\nk}_1(t)\leq m_\nk$ and $T_{gel}$ is a bounded random
variable by Proposition~\ref{gege}, we easily deduce from (\ref{1i2})
that $E \left[ \int_0^{T_{gel}} (M^{\nk}_1(t)/m_\nk) \ dt
\right]\longrightarrow 0$ as $k\to\infty$. Thus (\ref{1i4}) follows
from (\ref{chch}). 

\medskip

We next establish (\ref{1i1}) and (\ref{1i3}) and split the proof into
five steps. In the first two steps we show that, for $t>T_{gel}$, at
least one particle has a size of order $\delta m_n$ for some
$\delta>0$. Since such a particle is very attractive, we deduce in
Step~3 that no other large particle can exist and obtain
(\ref{1i3}). We then conclude in the last two steps that, for
$t>T_{gel}$, this single giant particle is solely responsible for the
loss of mass and obtain (\ref{1i1}). 

\smallskip

{\it Step~1.} Let $(\alpha_n)_{n\ge 1}$ be any sequence of positive
numbers such that $\alpha_n/m_n\to 0$ as $n\to\infty$. The aim of this
step is to show that 
\begin{equation}
\label{lanfeust}
\lim_{k\to\infty} P[M_1^{n_k}(T_{gel}+\e)>\alpha_{n_k} ] =1 \;\;\mbox{
for all }\;\; \e>0\,. 
\end{equation}
For that purpose, we introduce the stopping time
$$
\tau_k := \inf \left\{ t\geq 0 \; : \;\; M_1^{\nk}(t)> \alpha_\nk
\right\}\,, 
$$
and notice that we may assume that $\alpha_n\to\infty$ as $n\to\infty$
without loss of generality. Owing to the time monotonicity of
$M_1^{n}$, it suffices to prove that $P[\tau_{k}\geq
T_{gel}+\e]\longrightarrow 0$ as $k\to\infty$ for all $\e>0$ to
establish (\ref{lanfeust}). Assume thus for contradiction that there
is $\e>0$ such that $\delta:=\limsup_k P[\tau_{k}\geq
T_{gel}+\e]>0$. Then, on the one hand, we have $P[\tilde \Omega] \geq
\delta>0$ with $\tilde \Omega := \{\limsup_k \tau_{k} \geq
T_{gel}+\e\}$. On the other hand, for each $k\geq 1$, it is clearly
possible to couple the Marcus-Lushnikov processes
$(\mu^{n_k,m_{n_k}}_t)_{t\geq 0}$ and
$(\mu^{n_{k},\alpha_{n_{k}}}_t)_{t\geq 0}$ in such a way that they
coincide on $[0,\tau_{k})$. Hence, a.s. on $\tilde \Omega$, up to
extraction of a subsequence, 
$$
\lim_{k\to\infty} \sup_{[0,T_{gel}+\e/2]}
d(\mu^{n_{k},m_{n_{k}}}_t,\mu_t)= 
\lim_{k\to\infty} \sup_{[0,T_{gel}+\e/2]}
d(\mu^{n_{k},\alpha_{n_{k}}}_t,\mu_t)=0. 
$$
By Theorem~\ref{result} Points~1-(i) and~2, we deduce that the limit
$(\mu_t)_{t\geq 0}$ solves simultaneously the Flory and Smoluchowski
equations on $[0,T_{gel}+\e/2)$ with positive probability, which
contradicts Corollary~\ref{flnesm}. 

\smallskip  

{\it Step~2.} We now deduce from Step 1 that
\begin{equation}
\label{cixi}
\lim_{\delta\to 0} \liminf_k P\left[M_1^{n_k}(T_{gel}+\e)>\delta
m_{n_k} \right]=1 \;\;\mbox{ for all }\;\; \e>0\,. 
\end{equation} 
Assume for contradiction that there is $\e>0$ for which (\ref{cixi})
fails to be true. Then there exists $\gamma \in [0,1)$ such that
$\liminf_k P\left[M_1^{n_k}(T_{gel}+\e)>\delta m_{n_k}  \right]<
\gamma$ for all $\delta>0$. We may thus find a strictly increasing
sequence $(k_l)_{l\geq 1}$ such that
$P\left[M_1^{n_{k_l}}(T_{gel}+\e)> m_{n_{k_l}}/l \right]\le\gamma$ for
every $l \geq 1$. We then put $\alpha_{n_{k_l}}=m_{n_{k_l}}/l$ for
$l\ge 1$ (and e.g. $\alpha_n=m_n^{1/2}$ if $n\not\in\{ n_{k_l} : l\ge
1\}$). Then $\alpha_n/m_n\to 0$ as $n\to\infty$ and the assertion
(\ref{lanfeust}) established in Step~1  warrants that
$P\left[M_1^{n_k}(T_{gel}+\e)>\alpha_{n_k} \right]\longrightarrow 1$
as $k\to\infty$. But
$P\left[M_1^{n_{k_l}}(T_{gel}+\e)>\alpha_{n_{k_l}} \right]\leq
\gamma<1$ for all $l\ge 1$, hence a contradiction. 

\smallskip

{\it Step 3.} We are now in a position to prove (\ref{1i3}) which
somehow means that the other particles are small in the sense that 
\begin{equation}\label{cqvsp}
\lim_{b \to \infty} \limsup_{k\to\infty} P \left[ \left(
\int_{T_{gel}}^\infty X^{n_k}(s,b)\ ds \right) \geq \beta \right] =0 
\end{equation}
for all $\beta>0$ with the notation
$$
X^n(s,b) := \frac{1}{m_n} \sum_{i\geq 2} M^{n}_i(s)
\indiq_{[b,\infty)}\left( M^{n}_i(s) \right)\,. 
$$
First note that a.s., for all $s\geq 0$ and $n\geq 1$, we have
$$
\frac{1}{m_n^2} \sum_{i\ne j} M_i^{n}(s) M_j^{n}(s)
\indiq_{[b,m_n]}\left( M_i^{n}(s) \right) \indiq_{[b,m_n]}\left(
M_j^{n}(s) \right) \geq \frac{M_1^{n}(s)}{m_n} \indiq_{[b,\infty)}
\left( M_1^{n}(s) \right)\ X^n(s,b)\,, 
$$
since $M_i^{n}(s)\leq m_n$ for all $s\geq 0$ and $i\ge 1$. By
Lemma~\ref{fund} (with $a=m_n$), we obtain 
\begin{equation}\label{qqq}
E \left[ \int_{T_{gel}}^\infty \frac{M_1^{n}(s)}{m_n}
\indiq_{[b,\infty)}\left( M_1^{n}(s) \right)\ X^n(s,b)\ ds \right]\leq
\frac{L}{b^{\alpha}}. 
\end{equation}
 
We next fix $\beta>0$, $\eta>0$, and $b>0$. By (\ref{cixi}) there is
$\delta>0$ such that $\liminf_k P[M_1^{n_k}(T_{gel}+\beta/2) \geq
\delta m_{n_k}]\geq 1-\eta$. Recalling that $t\mapsto M_1^{n}(t)$ is
a.s. non-decreasing and $X^n(s,b)\leq 1$ for all $s\ge 0$ a.s., we
have for $k$ sufficiently large such that $\delta m_{n_k}>b$, 
\begin{eqnarray*}
&&P\left[\left(\int_{T_{gel}}^\infty X^\nk(s,b)\ ds \right) \geq \beta
\right]\leq P\left[\left(\int_{T_{gel}+\beta/2}^\infty X^\nk(s,b)\ ds
\right) \geq \beta/2 \right] \\ 
&&\leq P\left[ M_1^{n_k}(T_{gel}+\beta/2) \le \delta m_{n_k} \right] \\
&&+ P\left[ \left( \int_{T_{gel}+\beta/2}^\infty
\frac{M_1^{n_k}(T_{gel}+\beta/2)}{\delta m_{n_k}}
\indiq_{[b,\infty)}\left( M_1^{\nk}(T_{gel}+\beta/2) \right)
X^{n_k}(s,b)\ ds \right) \ge \beta/2 \right] \\ 
&&\leq 1 - P\left[ M_1^{n_k}(T_{gel}+\beta/2) > \delta m_{n_k} \right]
\\  
&& + P\left[ \left( \int_{T_{gel}+\beta/2}^\infty
\frac{M_1^{n_k}(s)}{m_{n_k}} \indiq_{[b,\infty)}\left( M_1^{\nk}(s)
\right) X^{n_k}(s,b)\ ds \right) \ge \beta\delta/2 \right]\\ 
&&\leq 1 - P\left[ M_1^{n_k}(T_{gel}+\beta/2) > \delta m_{n_k} \right]
+ \frac{2L}{b^\alpha\beta\delta}\,, 
\end{eqnarray*}
the last inequality being a consequence of (\ref{qqq}). Letting
$k\to\infty$ in the above inequality, we obtain, thanks to the choice
of $\delta$, 
$$
\limsup_{k\to\infty} P\left[\left(\int_{T_{gel}}^\infty X^\nk(s,b)\ ds
\right) \geq \beta \right] \le \eta +
\frac{2L}{b^\alpha\beta\delta}\,. 
$$
Now, we first pass to the limit as $b\to\infty$ and then as $\eta\to
0$ in the above inequality to obtain (\ref{cqvsp}), i.e. (\ref{1i3}). 

\smallskip

{\it Step 4.} Set $\gamma_t:=1 -\lc \mu_{T_{gel}+t}(dx),x\rc$ and
$B_{k}(t):=M^{\nk}_1(T_{gel}+t)/m_\nk$ for $t\ge 0$ and $k\ge 1$. Our
aim in this step is to prove that 
\begin{equation}\label{ccqqvv}
\lim_{k\to\infty} E\left[ \int_0^T |B_k(t) - \gamma_t|\ dt \right]=0
\;\;\mbox{ for all }\;\; T>0\,. 
\end{equation}
As before let $(\xhi_b)_{b>0}$ be a family of continuous
non-increasing functions such that $\xhi_b(x)=1$ for $x\in [0,b]$ and
$\xhi_b(x)=0$ for $x\geq b+1$. We then put 
$$
\gamma_t^b := 1-\lc \mu_{T_{gel}+t}(dx),x\xhi_b(x)\rc \;\;\mbox{ and
}\;\; \gamma_t^{b,k}:=1-\lc \mu^{\nk,m_\nk}_{T_{gel}+t}(dx),x\xhi_b(x)\rc
$$
for $b>0$, $k\ge 1$, and $t\ge 0$. On the one hand, we have a.s. that
$\gamma_t^b \longrightarrow \gamma_t$ as $b\to\infty$ for all $t\geq
0$. Since $|\gamma_t|+|\gamma_t^b|\leq 2$, we deduce from the Lebesgue
dominated convergence theorem that   
\begin{equation}
\label{donald}
\lim_{b \to \infty} E\left[\int_0^T |\gamma_t - \gamma_t^b|\ dt
\right] =0\,. 
\end{equation}
On the other hand, owing to the compactness of the support of
$\xhi_b$, we infer from (\ref{marsu}) that
$\gamma_t^{b,k}\longrightarrow \gamma_t^b$ a.s. for all $b>0$ and
$t\ge 0$. As $|\gamma_t^b|+|\gamma_t^{b,k}| \leq 2$, we use again the
Lebesgue dominated convergence theorem to obtain that 
\begin{equation}
\label{mickey}
\lim_{k \to \infty} E\left[ \int_0^T |\gamma_t^b - \gamma_t^{b,k}|\ dt\right] 
=0 \;\;\mbox{ for each }\;\; b>0\,.
\end{equation}
But $\gamma_t^{b,k}=A_{b,k}(t) + B_{k}(t)-C_{b,k}(t)$ a.s., where
\begin{eqnarray*}
A_{b,k}(t)&:=&\frac{1}{m_\nk}\sum_{i\geq 2} M^{\nk,m_\nk}_i(T_{gel}+t)
\left( 1-\xhi_b( M^{\nk,m_\nk}_i(T_{gel}+t))\right), \\
C_{b,k}(t) &:=& \frac{M^{\nk,m_\nk}_1(T_{gel}+t)}{m_\nk}
\xhi_b( M^{\nk,m_\nk}_1(T_{gel}+t)) \leq \frac{b+1}{m_\nk}.
\end{eqnarray*}
Clearly, 
\begin{equation}
\label{dingo}
\lim_{k\to\infty} E\left[\int_0^T C_{b,k}(t)\ dt \right]=0 \;\;\mbox{
for each }\;\; b>0\,,
\end{equation}
while, since $0\leq A_{b,k}(t)\leq X^\nk(T_{gel}+t,b)\leq 1$ a.s.,
\begin{equation}\label{prov}
\lim_{b\to \infty} \limsup_{k\to\infty} E\left[\int_0^T A_{b,k}(t)\
dt\right] = 0
\end{equation}
by (\ref{cqvsp}) and the Lebesgue dominated convergence theorem. 

Now, since $B_k(t) - \gamma_t= C_{b,k}(t) - A_{b,k}(t) +
(\gamma^{b,k}_t -\gamma^b_t) 
+ (\gamma^{b}_t -\gamma_t)$ for $b>0$, $k\ge 1$, and $t\ge 0$, it
follows from (\ref{mickey}) and (\ref{dingo}) that 
$$ 
\limsup_{k\to\infty} E\left[\int_0^T |B_k(t) - \gamma_t|\ dt
\right] \leq \limsup_{k\to\infty} E\left[\int_0^T A_{b,k}(t)\ dt
\right] + E\left[\int_0^T |\gamma_t - \gamma_t^b|\ dt\right]
$$
for all $b>0$. Letting $b\to \infty$ and using (\ref{donald}) and
(\ref{prov}) give (\ref{ccqqvv}). 

\smallskip

{\it Step 5.} To complete the proof of (\ref{1i1}), it remains to show
that, for all $\e>0$ and $\eta>0$, we have
\begin{equation}\label{finalhere}
\lim_{k\to\infty} P\left[\sup_{t\in [\eta,\infty)} |B_k(t) - \gamma_t|
\geq \e \right]=0. 
\end{equation}
Indeed, (\ref{finalhere}) clearly implies (\ref{1i1}) by the Lebesgue
dominated convergence theorem since $B_k(t)\le 1$ and $\gamma_t\leq 1$
for $t\ge 0$ a.s.\\
We thus fix $\e>0$ and $\eta>0$. Since $(\mu_t)_{t\ge 0}$ solves the
Flory equation (F) by Theorem~\ref{result} Point~1-(i), it follows
from Proposition~\ref{gege} that $t\longmapsto \gamma_t$ is
a.s. increasing and continuous on $[\eta,\infty)$ and
$\gamma_t\longrightarrow 1$ as $t\to\infty$ a.s. It is also
straightforward to check that $t\mapsto B_k(t)$ is a.s. non-decreasing
on $[\eta,\infty)$ with $B_k(t)\longrightarrow 1$ as $t\to\infty$
a.s. As a consequence of the a.s. monotonicity and boundedness of
$t\longmapsto \gamma_t$ and $B_k(t)$, we have for $t\ge T$
\begin{eqnarray*}
|B_k(t)-\gamma_t| & = & \max{\{ B_k(t) - \gamma_T + \gamma_T -
\gamma_t , \gamma_t - \gamma_T + \gamma_T - B_k(T) + B_k(T) - B_k(t)
\}} \\
& \le & \max{\{ 1 - \gamma_T , 1 - \gamma_T + \gamma_T - B_k(T) \}} \\
& \le & 1 - \gamma_T + \max{\{ 0 , \gamma_T - B_k(T) \}}\,,
\end{eqnarray*}
hence
\begin{equation}\label{jab3}
\sup_{t\in [T,\infty)}|B_k(t)-\gamma_t| \leq 1-\gamma_T +
|B_k(T)-\gamma_T| \;\;\mbox{ for all }\;\; T>0\,.
\end{equation}

To go further we will use the following result which resembles Dini's
theorem.
\begin{lemma}\label{diniL1}
Let $T>0$ and $f\in C([0,T])$ be a non-decreasing function. If
$(f_k)_{k\ge 1}$ is a sequence of non-decreasing functions on
$(0,T)$ such that $f_k\longrightarrow f$ in $L^1(0,T)$ as
$k\to\infty$, then $f_k\longrightarrow f$ in $L^\infty(\delta,T-\delta)$ as
$k\to\infty$ for every $\delta\in (0,T/2)$. 
\end{lemma}

Let $T>\eta$. By (\ref{ccqqvv}) and Proposition~\ref{gege},
$(B_k)_{k\ge 1}$ is a sequence of non-decreasing functions that
converges to the continuous and non-decreasing function $t\mapsto
\gamma_t$ in $L^1(0,T+\eta)$ a.s. and we use Lemma~\ref{diniL1} to
conclude that 
\begin{equation}\label{jab1}
\lim_{k\to\infty} P \left[\sup_{t\in [\eta,T]} |B_k(t)-\gamma_t| \geq
\e/2 \right] =0. 
\end{equation}  
We now infer from (\ref{jab3}) and (\ref{jab1}) that 
\begin{eqnarray*}
P\left[\sup_{t\in [\eta,\infty)}|B_k(t)-\gamma_t|\geq \e \right]
&\leq & P\left[\sup_{t\in [\eta,T]}|B_k(t)-\gamma_t| +
|1-\gamma_T|\geq \e \right]\\ 
&\leq & P\left[ 1 - \gamma_T \ge \e/2 \right] +
P\left[\sup_{t\in [\eta,T]}|B_k(t)-\gamma_t|\geq \e/2\right]\,,\\ 
\limsup_{k\to\infty} P\left[\sup_{t\in
[\eta,\infty)}|B_k(t)-\gamma_t|\geq \e \right] & \leq & P\left[ 1 -
\gamma_T \ge \e/2 \right]\,. 
\end{eqnarray*}
The above inequality being valid for any $T>\eta$, we may let
$T\to\infty$ to deduce (\ref{finalhere}) since $\gamma_T\longrightarrow
1$ as $T\to\infty$ a.s. 
\end{proof}

\medskip

We finally turn to the proof of the last statement of
Theorem~\ref{result}.

\medskip 

\begin{proof} {\bf of Point~3 of Theorem~\ref{result}.}
Here $(a_n)_{n\geq 1}$ is a sequence of positive real numbers such
that $a_n\to\infty$ and $a_n/m_n\to\gamma\in (0,1)$ as
$n\to\infty$. We consider the limit $(\mu_t)_{t\geq 0}$ of a
subsequence $\left\{ (\mu^{n_k,a_{n_k}}_t)_{t\geq 0} \right\}_{k\ge
1}$ in the sense that a.s. 
\begin{equation}
\label{lebrac}
\lim_{k\to\infty} \sup_{[0,T]} d(\mu^{n_k,a_{n_k}}_t,\mu_t)=0
\;\;\mbox{ for all }\;\; T\geq 0\,, 
\end{equation}
the existence of such a limit following by the same arguments 
as (\ref{prunelle}).

We then introduce 
$$
T_1 := \inf\left\{ t\geq 0\; : \; \lc \mu_0(dx) - \mu_t(dx),x \rc \geq
\gamma\right\}\,, 
$$
and aim at showing that $(\mu_{t})_{t\in [0, T_1)}$ solves a.s. the
Flory equation $(F)$.

For $n\ge 1$, we set
$$
T_1^n := \inf \{ t \geq 0\; ; \; \lc \mu^{n,a_n}_t,
\indiq_{(a_n,\infty)} \rc >0 \}\,, 
$$
which represents the first time that a particle of size exceeding
$a_n$ appears in the Marcus-Lushnikov process $(\mu^{n,a_n}_t)_{t\geq
0}$. For each $n\ge 1$, it is clearly possible to build a classical
Marcus-Lushnikov process $(\mu^{n,m_n}_t)_{t\geq 0}$ (i.e. without
cut-off) such that $\mu^{n,m_n}_t=\mu^{n,a_n}_t$ for $t\in[0,T_1^n]$
a.s. In particular we have also $T_1^n = \inf \{ t \geq 0\; : \; \lc
\mu^{n,m_n}_t(dx), \indiq_{(a_n,\infty)} \rc >0 \}$ a.s. Denoting by
$M_1^n(t)$ the size of the largest particle at time $t$ in the process
$(\mu^{n,m_n}_t)_{t\geq 0}$, we clearly have $T_1^n = \inf \{ t \geq
0\; : \; M_1^n(t)>a_n\}$ a.s.

By the tightness result of Proposition~\ref{tension}, we may assume
that, up to extracting a further subsequence (not relabeled),
$(\mu^{n_k,m_\nk}_t)_{t\geq 0}$ converges to $(\nu_t)_{t\geq 0}$ in 
$\mathbb{D}([0,\infty), \mes)$. By the Skorokhod representation
theorem and Lemma~\ref{toujoursok}, we may assume that a.s.
\begin{equation}
\label{clarabelle}
\lim_{k\to\infty} \sup_{[0,T]} d(\mu^{n_k,m_\nk}_t,\nu_t)=0 \;\;\mbox{
for all }\;\; T>0\,.
\end{equation}

By Theorem~\ref{result} Point~1-(i), $(\nu_t)_{t\geq 0}$ solves a.s.
the Flory equation $(F)$. Introducing 
$$
S_1 := \inf\left\{ t\geq 0\; : \; \lc \mu_0(dx) - \nu_t(dx),x \rc 
\geq \gamma\right\}\,,
$$ 
we claim that  
\begin{equation}\label{cqqv}
\lim_{k\to\infty} P[|S_1 - T_1^\nk|> \e]=0 \;\;\mbox{ for 
all }\;\; \e>0\,.
\end{equation}

Taking (\ref{cqqv}) for granted, we deduce that $\mu_t=\nu_t$ for
$t\in [0,S_1)$ a.s. since $\mu^{n,m_n}_t=\mu^{n,a_n}_t$ for
$t\in[0,T_1^n]$ a.s. for all $n\geq 1$. This implies that $S_1=T_1$ 
a.s., because the subset $\{\pi \in \mes; \; \lc\mu_0(dx) - 
\pi(dx),x  \rc \geq \gamma\}=\{\pi \in \mes; \; \lc \pi(dx),x  
\rc \leq 1- \gamma\}$ is closed in $\mes$ endowed with 
the vague topology, and both $t\mapsto \mu_t$ and $t\mapsto \nu_t$ 
are a.s. continuous for that topology by Lemma~\ref{toujoursok}.
Therefore, $(\mu_t)_{t\in[0,T_1)}$ solves a.s. the Flory equation.

\smallskip 

We are left with the proof of (\ref{cqqv}). To this end we will use
(\ref{1i1}) and (\ref{1i2}) (with the weak limit $(\nu_t)_{t\geq 0}$ 
of the classical Marcus-Lushnikov process $\left( \mu_t^{n_k,m_{n_k}} 
\right)_{k\ge 1}$). Introducing the (random) gelling time $S_{gel}$ 
of $(\nu_t)_{t\ge 0}$ given by
$$
S_{gel} := \inf\{t\geq 0\; : \; \lc \nu_t(dx),x\rc < 
\lc \mu_0(dx),x\rc\}\,,
$$
we recall that a.s. the map $t \longmapsto \lc \nu_t(dx),x\rc$ is 
constant and equal to $1$ on $[0,S_{gel})$ and continuous and 
decreasing on $(S_{gel},\infty)$ by Proposition~\ref{gege}. In the 
proof of (\ref{cqqv}), we have to handle separately the events 
$S_{gel}<S_1$ and $S_{gel}=S_1$, the latter being not ruled out 
{\it a priori} due to the possible discontinuity of $t\longmapsto 
\lc \nu_t(dx),x\rc$ at $t=S_{gel}$. 

Fix $\e>0$ and write 
$$
P[|S_1 - T_1^\nk|>\e] = P[U_k]+P[V_k]+P[W_k]
$$
with
\begin{eqnarray*}
U_k & := & \left\{S_{gel} \le S_1 \le S_{gel}+\e/2\;,\;\; T_1^\nk 
< S_1-\e \right\}\,,\\
V_k & := & \left\{S_{gel} +\e/2 < S_1\;,\;\; T_1^\nk < 
S_1-\e \right\}\,,\\
W_k & := & \left\{S_{gel} \le S_1\;,\;\; T_1^\nk > S_{1}+\e 
\right\}\,.
\end{eqnarray*}

First, on $U_k$, we have $a_{n_k} \leq M_1^{n_k}(S_1-\e) \leq M_1^{n_k}(S_{gel}-\e/2)$, so that 
$$
P[U_k]\le P\left[ M_1^{n_k}(S_{gel}-\e/2) \geq a_{n_k}\right] 
\leq P\left[ \frac{M_1^{n_k}(S_{gel}-\e/2)}{m_{n_k}} \geq 
\frac{a_{n_k}}{m_{n_k}}\right] 
\mathop{\longrightarrow}_{k\to\infty} 0
$$ 
by (\ref{1i2}) since $a_n/m_n\to \gamma>0$ as $n\to\infty$.

Next, introducing $\tau:= S_1 -\e/4$ and $Z:=\lc \nu_\tau(dx),x\rc - 1 + \gamma$, it follows from the a.s. strict monotonicity of 
$t \longmapsto \lc \nu_t(dx),x\rc$ on $(S_{gel},\infty)$
and the definitions of $T_1^n$ and $S_1$ that 
$$
Z>0 \;\;\mbox{ a.s. and }\;\; V_k\subset \left\{ \frac{M_1^{n_k}(\tau)}{m_{n_k}} - \left( 1 - \lc \nu_\tau(dx) , x \rc \right) \ge \frac{a_{n_k}}{m_{n_k}} - \gamma + Z \right\}\,.
$$
Let $\eta>0$. For $k$ large enough we have $|\gamma - a_{n_k}/m_{n_k}|\le \eta/2$ and since $\tau>S_{gel}+\e/4$ a.s. \begin{eqnarray*}
& & E\left[ \sup_{t\in [S_{gel}+\e/4,\infty)}{\left\{ \frac{M_1^{n_k}(t)}{m_{n_k}} - \left( 1 - \lc \nu_t(dx) , x \rc \right) \right\}} \right] \\
& \ge & E\left[ \indiq_{V_k}\ \left( \frac{M_1^{n_k}(\tau)}{m_{n_k}} - \left( 1 - \lc \nu_\tau(dx) , x \rc \right) \right) \right] \ge E\left[ \indiq_{V_k}\ \left( Z + \frac{a_{n_k}}{m_{n_k}} - \gamma \right) \right] 
\geq E \left[\indiq_{V_k} (Z-\eta/2) \right]\\
& \ge & \frac{\eta}{2}\ E\left[ \indiq_{V_k}\ \indiq_{[\eta,\infty)}(Z) \right] - \frac{\eta}{2}\ E\left[ \indiq_{(0,\eta)}(Z) \right] \ge \frac{\eta}{2}\ P[V_k] - \eta P[ Z\in (0,\eta) ]\,.
\end{eqnarray*}
Multiplying the above inequality by $2/\eta$ and letting $k\to\infty$ with the help of (\ref{1i1}) give
$$
\limsup_{k\to\infty} P[V_k] \le 2 P[ Z\in (0,\eta) ] \;\;\mbox{ for all } \eta>0\,.
$$
As $Z>0$ a.s., the right hand side of the above inequality converges to zero as $\eta\to 0$. Consequently, $P[V_k]\longrightarrow 0$ as $k\to\infty$. 

Similarly, introducing $\sigma:= S_1+\e$ and $Y:=1 -\lc \nu_\sigma(dx),x\rc
- \gamma$, the a.s. strict monotonicity of $t \longmapsto \lc \nu_t(dx),x\rc$  on $(S_{gel},\infty)$ and the definitions of $T_1^n$ and $S_1$ warrant that 
$$
Y>0 \;\;\mbox{ a.s. and }\;\; W_k \subset \left\{ \left( 1 - \lc \nu_\sigma(dx) , x \rc \right) - \frac{M_1^{n_k}(\sigma)}{m_{n_k}} \ge \gamma - \frac{a_{n_k}}{m_{n_k}} + Y \right\}\,.
$$
Arguing as for $V_k$, we have for all $\eta>0$ and $k$ large enough
$$
E\left[ \sup_{t\in [S_{gel}+\e,\infty)}{\left\{ \frac{M_1^{n_k}(t)}{m_{n_k}} - \left( 1 - \lc \nu_t(dx) , x \rc \right) \right\}} \right] \ge \frac{\eta}{2}\ P[W_k] - \eta P[ Y\in (0,\eta) ]\,.
$$
We then proceed as before to deduce from (\ref{1i1}) and the a.s. positivity of $Y$ that $P[W_k]\longrightarrow 0$ as $k\to\infty$ and thus complete the proof of (\ref{cqqv}).
\end{proof} 

\section{Numerical illustrations}\label{numerical}\setcounter{equation}{0}

We consider the {\it monodisperse} initial condition $\mu_0=\delta_1$
and the {\it multiplicative kernel} $K(x,y)=xy$. Under these conditions, 
there is an explicit solution to the Smoluchowski equation $(S)$
given by 
$$
\hat \mu_t(dx) := \sum_{k\geq 1} \hat c(t,k) \delta_k(dx) \;\;\mbox{ with }\;\; \hat c(t,k) := \left\{
\begin{array}{lcl}
\displaystyle{\frac{k^{k-2}}{k!}t^{k-1}e^{-kt}} & \;\mbox{ for }\; & t\in [0,1]\,, \\
& & \\
\displaystyle{\frac{k^{k-2}}{k!}t^{-1}e^{-k}} & \;\mbox{ for }\; & t\geq 1\,.
\end{array}
\right.
$$
For the same initial condition, the Flory equation $(F)$ has also an 
explicit solution given by 
$$
\mu_t(dx) := \sum_{k\geq 1} c(t,k) \delta_k(dx)\;\;\mbox{ with }\;\; c(t,k) := \frac{k^{k-2}}{k!}t^{k-1}e^{-kt} \;\mbox{ for }\; t\geq 0\,.
$$
Before proceeding to simulations, let us point out that $\lc \mu_t(dx),x \rc = 1$ for $t\in [0,1]$, while $\lc \mu_t(dx),x \rc=t^*/t$ for $t > 1$, where $t^*\in (0,1)$ is the unique solution to $t^*e^{-t^*}=te^{-t}$ in $(0,1)$. Easy computations show that
$$
T_1(\gamma) :=\inf\{t\geq 0; \;\lc \mu_0(dx)-\mu_t(dx),x \rc \geq \gamma\} = -\frac{\ln (1-\gamma)}{\gamma} \;\;\mbox{ for }\;\; \gamma \in (0,1)\,.
$$

In Figures~1 to~5, the plain, dashed, and dotted lines represent $\mu^{n,a_n}_t(\{2\})$, $c(t,2)$, and $\hat c(t,2)$, respectively, as functions of $t$. We observe that, as explained by Theorem~\ref{result}, 
\begin{itemize}
\item[(i)] for $a_n \ll m_n$,  $\mu^{n,a_n}_t$ approximates the solution to the Smoluchowski equation, see Figure~1,
\item[(ii)] for $a_n=m_n$, $\mu^{n,a_n}_t$ approximates the solution to the Flory equation, see Figure~2,
\item[(iii)] for $a_n = \gamma m_n$ with $\gamma\in (0,1)$, $\mu^{n,a_n}_t$ approximates the solution to the Flory equation until the time $T_1(\gamma)$, and then changes its behaviour: see Figure~3 ($\gamma = 0.5$, $T_1(0.5)=1.386$), Figure~4 ($\gamma = 0.8$, $T_1(0.8)=2.012$) and Figure~5 ($\gamma = 0.33$, $T_1(0.33)= 1.21$). Note that Figure~5 shows that the behaviour of $\mu^{n,a_n}_t$ {\it bifurcates} at least twice on $t\in [0,3]$. The second bifurcation certainly corresponds to the time where a second {\it giant} particle with size $10^5$ appears.
\end{itemize}

\begin{figure}[l]
\centerline{Figure 1: $n=m_n=10^4$, $a_n=10^2$ 
\hskip1cm Figure 2: $n=m_n=10^4$, $a_n=10^4$}
\centerline{\includegraphics[height=3.8cm]{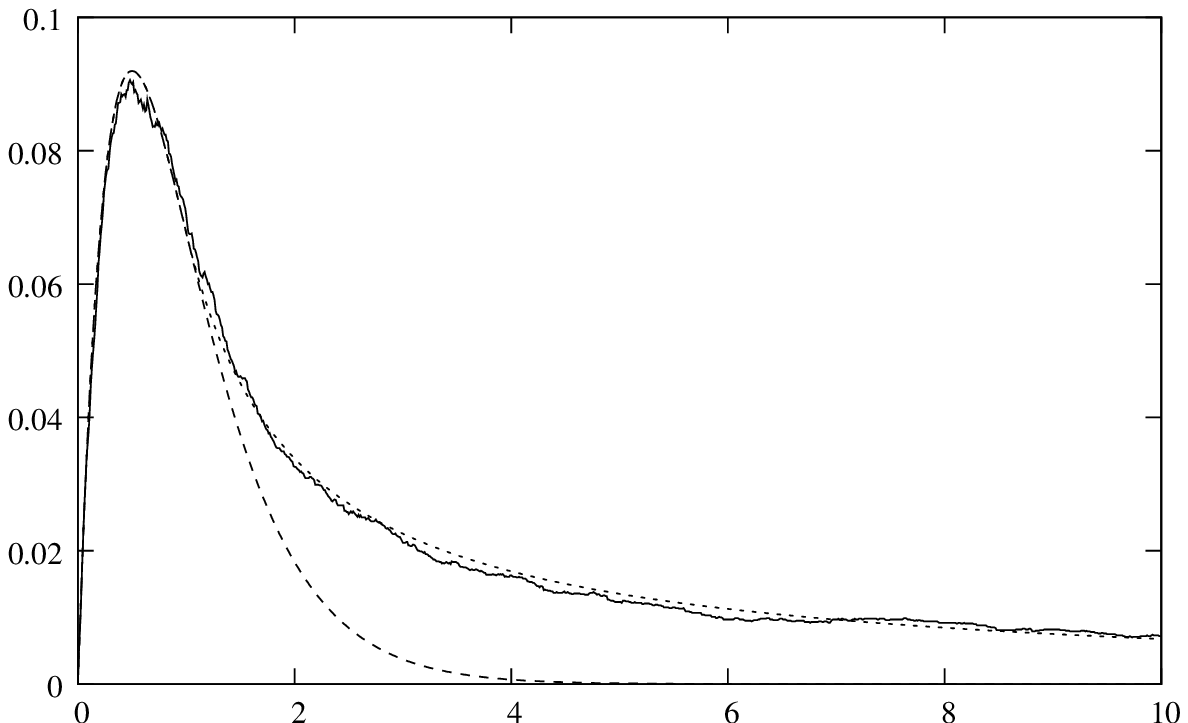} \hskip0.3cm
            \includegraphics[height=3.8cm]{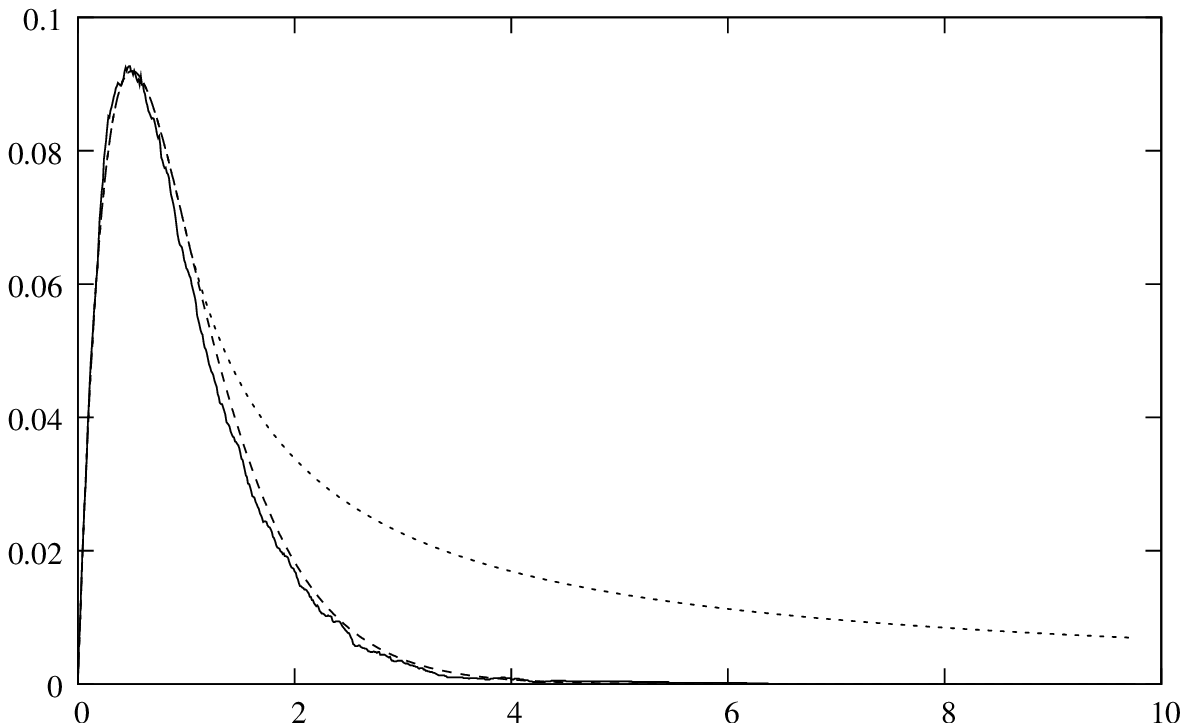}}
\end{figure}
\begin{figure}[l]
\centerline{Figure 3: $n=m_n=10^4$, $a_n=5.10^3$ \hskip1cm 
Figure 4:  $n=m_n=10^4$, $a_n=8.10^3$}
\centerline{\includegraphics[height=3.8cm]{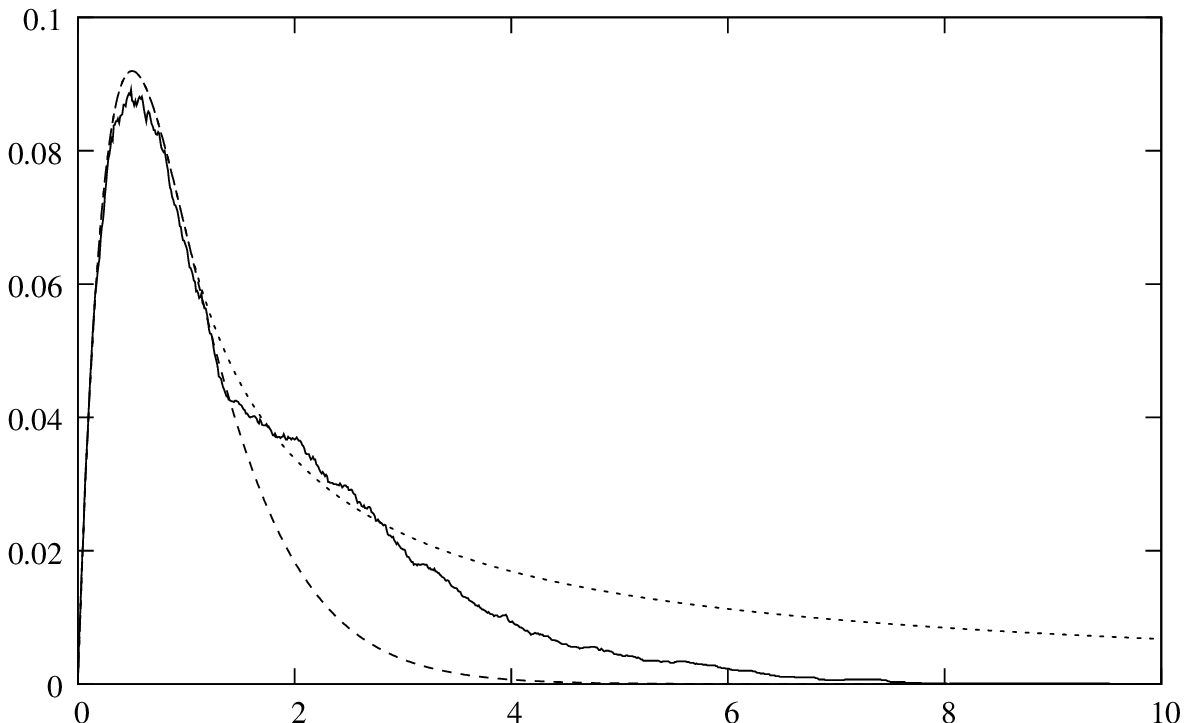} \hskip0.3cm
            \includegraphics[height=3.8cm]{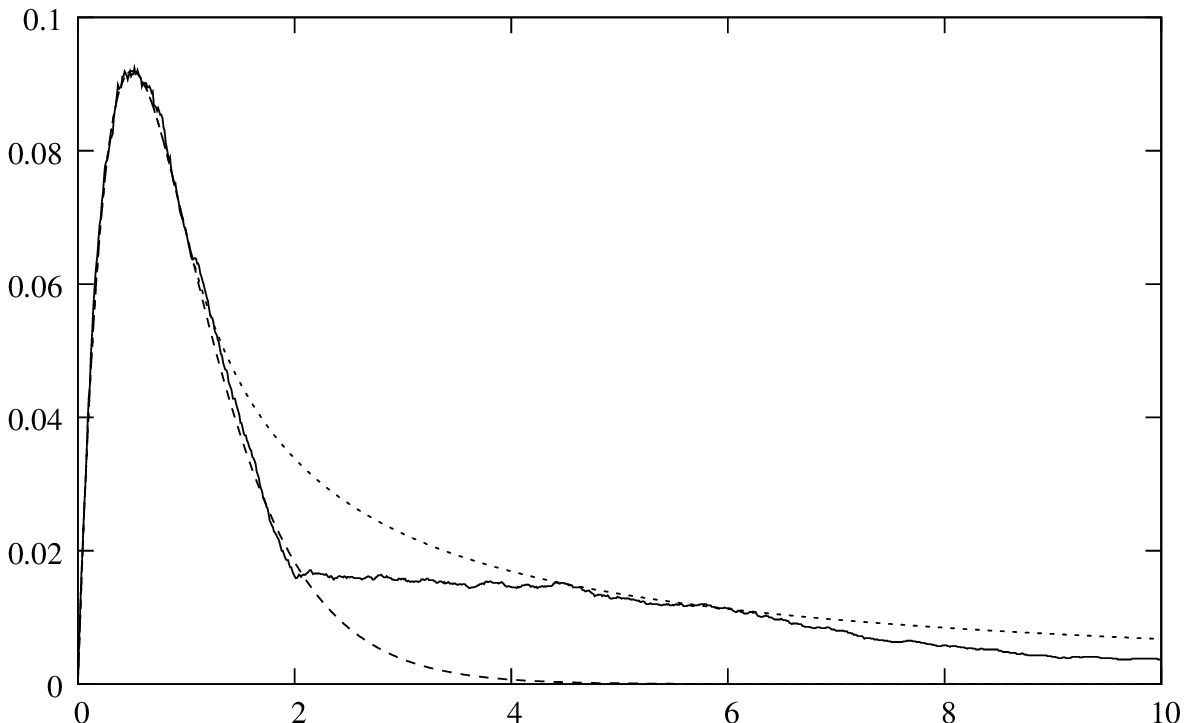}}
\end{figure}
\begin{figure}[l]
\centerline{Figure 5: $n=m_n=3.10^5$, $a_n=10^5$.}
\centerline{ \includegraphics[height=6cm]{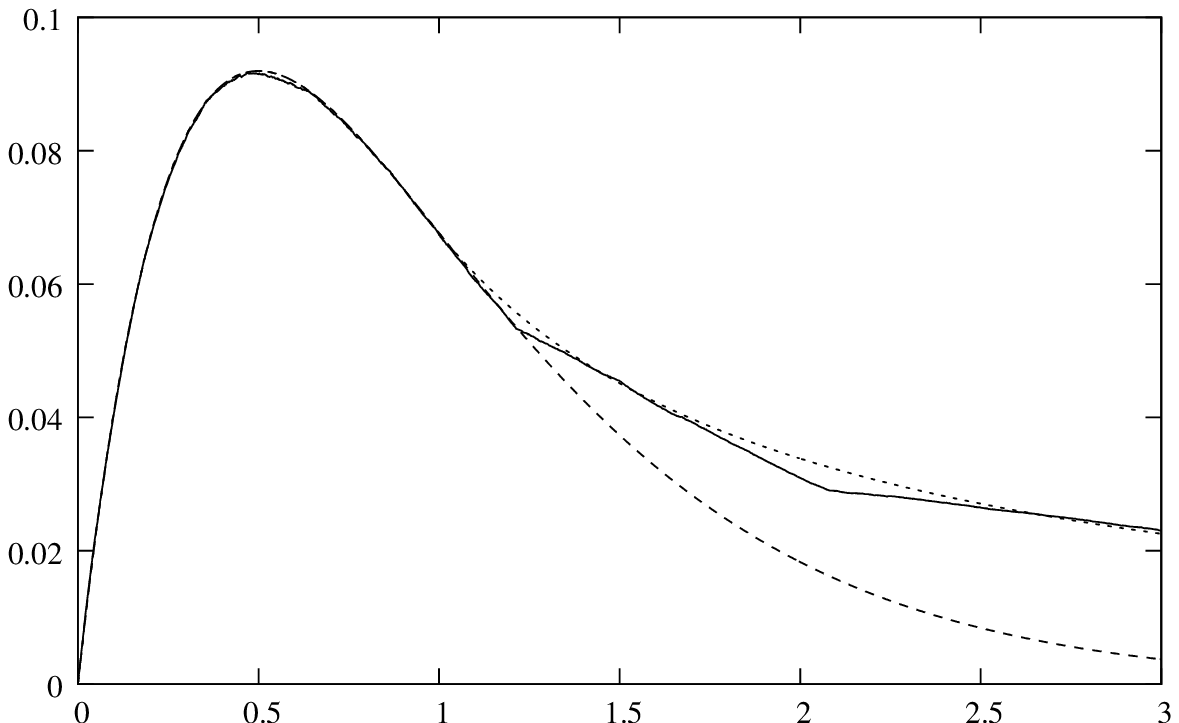} }
\end{figure}

\def\refname{References}

\end{document}